\documentclass[]{article} 
\usepackage{helvet}
\usepackage{url}              
\usepackage[utf8x]{inputenc}  
\usepackage[left=2cm,right=2cm,top=2cm,bottom=2cm]{geometry}  
\usepackage{xcolor}           
\usepackage{ulem}             
\usepackage{tikz}             
\usepackage{pdflscape}        
\usepackage{colortbl}         
\usepackage{graphicx}         
\usepackage{etoolbox}
\usepackage{overpic}
\usepackage{pict2e} 

\usepackage{amsmath}
\usepackage{amssymb}    
\usepackage[colorlinks=true, allcolors=blue]{hyperref}
\usepackage{caption}
\usepackage{subcaption}
\usepackage{float}
\usepackage{wrapfig}

\newcommand{\sign}{\text{sign}}

\patchcmd{\thebibliography}{\section*{\refname}}{}{}{}

\def\softd{{\leavevmode\setbox1=\hbox{d}%
          \hbox to 1.05\wd1{d\kern-0.4ex{\char039}\hss}}}

\title{Longitudinal wall shear stress evaluation using centerline projection approach in the numerical simulations of the patient-based carotid  artery}
\author{Kevin Richter$^1$\footnote{corresponding author, eMail: richter@uni-landau.de}, Tristan Probst$^1$,   Anna Hundertmark$^1$, \\ Pepe Eulzer$^2$, Kai Lawonn$^2$}
\date{%
    $^1$ Institute of Mathematics, Faculty of Natural and Enviromental Sciences,  University of Koblenz-Landau, Germany,\\
    $^2$  Faculty of Mathematics and Computer Science, University of Jena, Germany\\
    [2ex]%
    \today
}

\begin{document}

\maketitle

\begin{abstract}
In this numerical study areas  of the carotid bifurcation and  of a distal stenosis  in the internal carotid artery are closely observed to evaluate the patient's current risks of ischemic stroke. An indicator for the vessel wall defects is the stress the blood is exerting on the surrounding vessel tissue, expressed  standardly by the  amplitude of  the wall shear stress vector  (WSS) and its oscillatory shear index.
In contrast, our orientation-based shear evaluation detects negative shear stresses corresponding with  reversal flow appearing in low shear areas.  
In our investigations of longitudinal component of the wall shear vector, tangential vectors aligned longitudinally with the vessel are necessary. However, as a result of stenosed regions and  imaging segmentation techniques from patients' CTA scans,
the geometry model's mesh is  non-smooth on its surface areas and the automatically generated tangential vector field 
is discontinuous and multi-directional,  making an interpretation of the orientation-based risk indicators unreliable. 
We improve the evaluation of longitudinal shear stress  by applying the projection of the vessel's center-line to the surface to construct smooth tangetial field aligned longitudinaly with the vessel. We validate our approach for the longitudinal WSS component
and the corresponding oscillatory index by comparing them to results obtained using automatically generated tangents in both rigid and elastic vessel modeling as well as to amplitude based indicators.  The major benefit of our WSS evaluation based on  its  longitudinal component  for the cardiovascular risk assessment  is  the detection of negative WSS indicating persitent reversal flow. This is impossible in the case of  the amplitude-based WSS. 
\end{abstract}
\medskip    
\noindent
\\{\bf Keywords:}
hemodynamics, 
fluid-structure interaction, finite element method, cardiovascular risk indicators, 
longitudinal wall shear stress, oscillatory shear index.

\section{Introduction}
The importance of a healthy and functioning cardiovascular system is reflected in the WHO death statistics of 2019. Ischemic stroke was the disease responsible for the highest proportion of deaths across all countries and wealth levels \cite{who01}. The cause of ischemic stroke is an arterial vascular disease, which in its most common form, atherosclerosis, is an inflammatory response of the vessel wall to lipid metabolism disturbances and endothelial stress. This  leads to the formation of multi-focal plaques and thus to the narrowing and hardening of the arteries and consequently to an insufficient supply of oxygen to the brain \cite{debus}. A special role in atherosclerosis development plays the carotid artery, which is responsible for an estimated 18 - 25 \% of thromboembolic strokes \cite{iannuzzi}. In the carotid bifurcation the common carotid artery splits into the external and the internal carotid artery. While the former is responsible for supplying blood to the head and upper neck organs, the latter supplies blood to the brain. Both, the death toll of ischemic strokes and the drastic increase in the general prevalence of atherosclerosis, which is related to demographic change and the accompanying burden on health and care system, make it necessary to adequately address the danger posed by atherosclerosis. To provide necessary tools predicting the locations of sites susceptible to atherosclerotic damage, as well as to make recommendations for their optimal treatment, is one of the main goals of modern medicine.\\

\noindent
The predictions of atherosclerosis and further cardiovascular risk through numerical simulations have become popular over the last decades. Especially in the field of fluid-structure-interaction there is an inexhaustible range of publications. They span from the incorporation of different mathematical models, such as non-Newtonian fluids \cite{hundertmarkshear, janelashear}, different structure models, e.g. shell models and membrane models \cite{canicshell,ciarletshell}, to exploring new numerical methods, which effectively tackle the multi-physicality with splitting techniques \cite{bukacsplitting,Hundertmark2010,Hundertmark2013}, just to name a few.

For the purpose of risk quantification, parameters derived from numerical flow data are of great interest, e.g. wall shear stress (WSS) \cite{antiga, arzani, Giorgio_low_WSS, Giorgio_WSS2, taylor}, or oscillatory shear index (OSI) \cite{antiga,Qua_OSI, Tremmel_OSI}, both referred to be correlated with cardiovascular risk. 
It has been reported \cite{ku}, that apart from regions with high amplitude WSS, 
areas with low and temporary oscillating WSS promote atherosclerotic processes. The multi-directional behavior of WSS has also been linked to potential risk zones, see results on transversal WSS and other metrics in \cite{hoogendoorn}. 

Some of the recent results report on proper visualisation tools for the localisation of cardiovascular risk zones. These are based on the interplay of imaging techniques for exploring the vessel morphology and simulated data as velocity streamlines or WSS,  see, e.g., \cite{Eulzer1, Eulzer2} and further citations therein. For those tools and their underlying numerical simulations the patient's unique vessel morphology plays a crucial role for reliable risk predictions of the above mentioned parameters.

\medskip

The aim of our study is reliable numerical modeling and a quantification of the impact of the fluid flow dynamics on endothelial stresses in the carotid artery of a clinical  patient. The obtained  numerical data  
for a set of patients 
serve as the training set for machine learning algorithms within the research project MLgSA, see Acknowledgement, to explore potential stroke risks.
The domain of interest is the carotid bifurcation area with its separation of the common carotid artery into the internal carotid artery (ICA) and the external carotid artery (ECA).  
The 3D lumen of the carotid vessel tree including its stenosed region, see Fig. \ref{Fig:Computational Geometry}, as well the shape of the surrounding wall tissue,  Fig. \ref{Fig:thick wall}, have been reconstructed from the computer tomography angiography (CTA) data set of a clinical patient with the method described in \cite{Eulzer1}. 
This carotid vessel tree as well as its arterial wall shape have been imported 
as fluid and solid regions  
in the finite-element-based  software \textit{Comsol Multiphysics}. 
We perform numerical simulations based on the incompressible fluid flow model for both rigid walls
as well as compliant walls. For the latter 
the  fluid-structure interaction (FSI) model describing the interplay of the fluid and thick compliant walls is considered. The corresponding wall tissue mechanics is modeled by the deformation of a linear elastic material. 

The obtained numerical results are used to evaluate the hemodynamic risk parameters: the wall shear stress (WSS) and the  oscillatory shear index (OSI), measuring the temporal change of the wall shear stress direction over one cardiac cycle.  
Herein, we consider  the longitudinal component of  WSS and compare it to  the results for the amplitude of the WSS vector. Moreover, we evaluate the corresponding oscillatory indices using the temporal mean of the longitudinal WSS component, as well as the amplitude of the mean WSS vector. We give their hemodynamic interpretation and discuss the benefits of the orientation-based WSS evaluation in view of reversal flow detection.

For the calculation of the longitudinal WSS value the choice of proper tangential vectors  
is crucial. The analogous sensitivity  of the transversal WSS to the orientation of the tangential field has been addressed in \cite{gallo, hoogendoorn, mohamied16, mohamied14}.
Regarding the realistic geometries with spatially non-smooth surface topology, the mesh-based erratic tangential vector field  
leads to problematic, spatially discontinuous behavior of longitudinal (or transversal) WSS on the surface of the extracted geometry. 
 We  address this problem  and improve the evaluation of the longitudinal component of  WSS in complex realistic geometries. 
 Our approach is based on choosing tangential vectors obtained by the projection of the centerline of the vessel tree to the vessel surface, previously used by, e.g. Morbiducci et al. \cite{morbiducci15} or Arzani $\&$ Shadden \cite{arzani}.  
We  present numerical data and two wall parameters WSS and OSI derived from the projection method and from the automatically generated, mesh-based tangential vector fields. Further, we compare the values obtained for rigid as well as compliant carotid vessel walls in order to examine the importance of the compliance of wall tissue in considered mathematical model of the carotid artery.


\section{Mathematical modeling of the carotid flow}

\noindent

\begin{figure}[h!]
  \centering
  \begin{minipage}[b]{.45\linewidth} 
      \includegraphics[width=\linewidth]{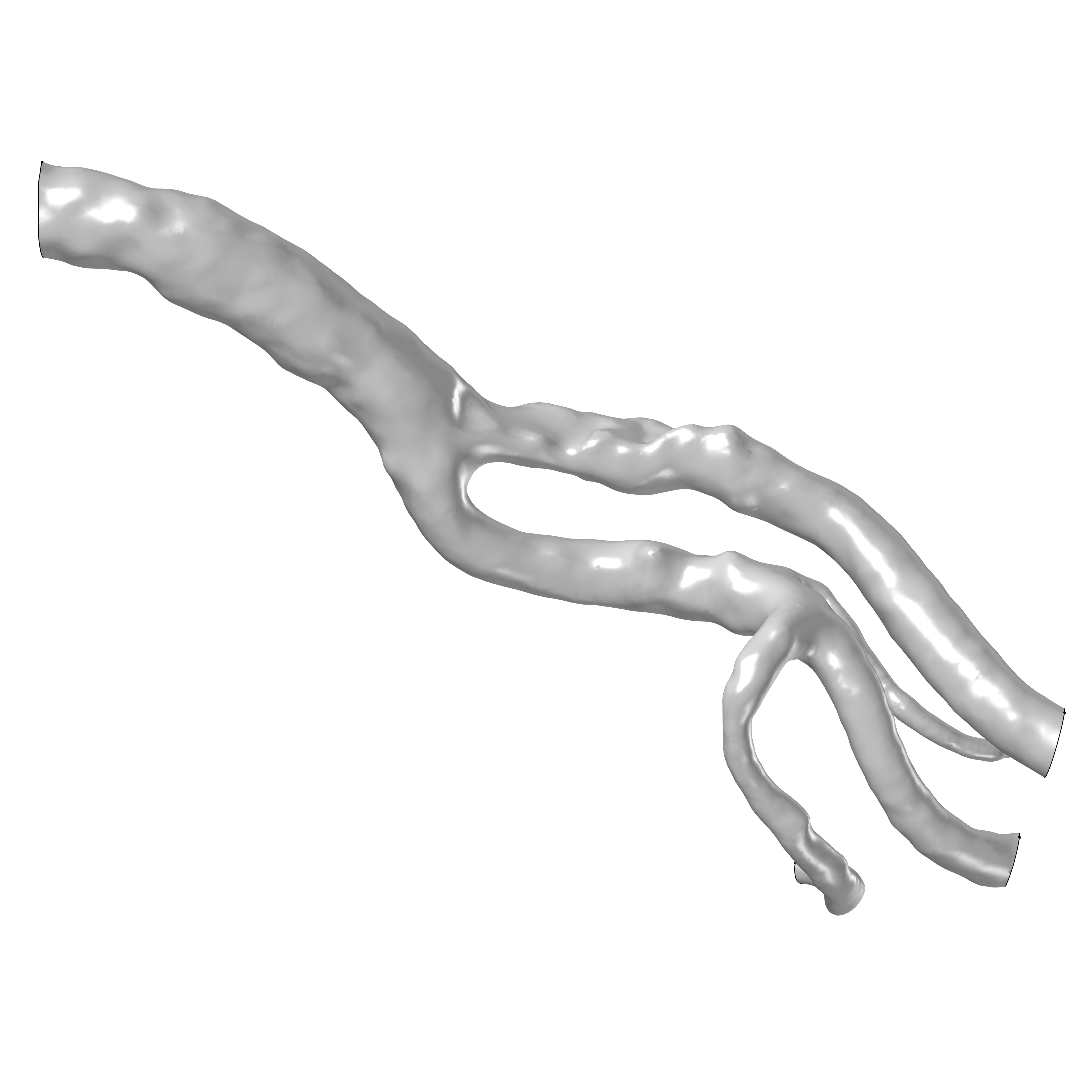}
      \caption{Computational geometry  $\Omega^f$ of the inner human carotid lumen for the fluid flow modeling with rigid walls.}
      \label{Fig:Computational Geometry}
  \end{minipage}
  \hspace{.05\linewidth}
  \begin{minipage}[b]{.45\linewidth} 
      \includegraphics[width=0.9\linewidth]{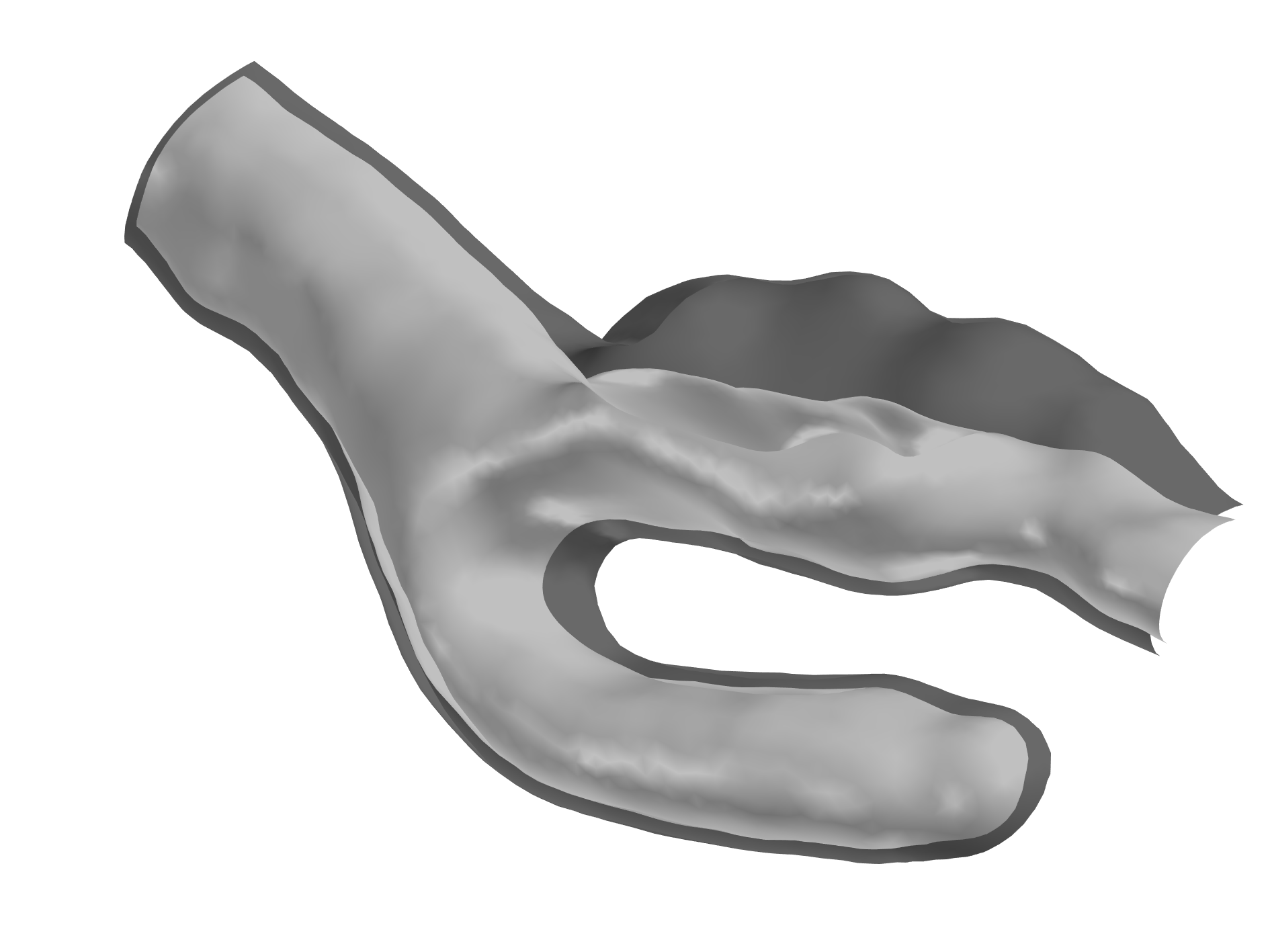}
      \caption{Computational geometry $\Omega^f_t \cup \Omega^s_t$ for the FSI-modeling of the human carotid including wall tissue, area of interest.}
      \label{Fig:thick wall}
  \end{minipage}
\end{figure}

To investigate the individual risk of arterial defects and imminent health issues the individual examination of the carotid geometry can give a more reliable assessment of the patient's current situation. A CTA scan of the patient is processed to a detailed three-dimensional arterial geometry, as described in \cite{Eulzer1} and the references therein. This technique was applied on the arterial lumen and expanded on the vessel wall domain as well.
The resulting
rigid domain of arterial lumen which includes part of the common carotid, its branching into internal and external carotid artery and two sub-branches of the latter, compare Fig.~\ref{Fig:Computational Geometry}, has been considered for numerical simulations. On the other hand, simulations which featured fluid-structure-interaction are supplemented with a vessel domain which also reveals stenotic regions and the distinct wall thickness distribution of the patient throughout the area of interest, see Fig.~\ref{Fig:thick wall}.
By cutting the computational domain proximal to the bifurcation area for the flow to fully develop and distal such that the boundary conditions don't affect the dynamics of the flow, the geometry is defined. 

\medskip

In particular, we denote with $\Omega_t=\Omega^f_t \cup \Omega^s_t\;,\;t\in[0,T]$ the deforming fluid and structural domain in time $t$, respectively and their shared boundary with $\Gamma^{fsi}_t=\Omega^f_t \cap \Omega^s_t$. 
%
%
To model the blood flow in a deforming vessel the incompressible Navier-Stokes equations in a moving domain are considered
in the arbitrary Lagrangian-Eulerian (ALE) formulation. They read as 
\begin{equation}\label{eq:ns}
    \rho_f\frac{D \boldsymbol{u}}{D t}+ \rho_f\left((\boldsymbol{u}-\boldsymbol{w}) \cdot  \nabla\right) \boldsymbol{u} 
    = \nabla \cdot \boldsymbol{T}_f, \quad \nabla \cdot \boldsymbol{u}=0 \qquad \text{in}\quad \Omega^f_t,
\end{equation}
with the Cauchy stress tensor  $\boldsymbol{T}_f=-p\boldsymbol{I}+2\mu\boldsymbol{D}(\boldsymbol{u})$ and the strain rate tensor $\boldsymbol{D}(\boldsymbol{u})=\frac1{2}(\nabla \boldsymbol{u}+ \nabla \boldsymbol{u} ^T)$. 
The fluid velocity $\boldsymbol{u}$ and pressure $p$ are solved in a fluid domain 
for constant viscosity $\mu$ and density $\rho_f$. 
Here $\boldsymbol{w}(X,t)=\frac{\partial x(X,t)}{\partial t},\  x\in  \Omega^f_t $ describes the fluid domain deformation velocity  with respect to the reference domain $\Omega^f_0$ and $\frac{D \boldsymbol{u}}{D t}$ is the total (material) time derivative of the fluid velocity.

The boundary of the fluid domain consists of inflow and outflow boundary as well as the shared fluid-structure interface,  $\Gamma^f=\Gamma^f_{in}\cup\Gamma^f_{out}\cup\Gamma^{fsi}_{t}$.
On the fixed inflow boundary $\Gamma^f_{in}$ a pulsating blood flow-rate was implemented, based on measured data taken from \cite{perktold95} as shown in Fig. \ref{fig:flow rate}. The outflow was assumed to have zero normal stress (do-nothing condition), that is
\begin{equation}
    \boldsymbol{T}_f\vec{n}_{out}=-p_{out}\vec{n}_{out}\equiv\boldsymbol{0}\qquad\text{ on }\Gamma_{out}^f.
\end{equation}

\noindent
The mechanics of the arterial wall is modeled by a linear elastic material which reads in the reference configuration as,
\begin{equation}
    \begin{split} \label{eq:def}
        \rho_s \frac{\partial^2 \boldsymbol{d}}{\partial t^2} &= \nabla \cdot (\boldsymbol{FS})^T+\boldsymbol{f}_s,\quad\text{ in } \Omega^s_0,\\
        \text{with}\quad
        \boldsymbol{S}&=J\boldsymbol{F}^{-1}\left(\boldsymbol{C:\epsilon}\right)\boldsymbol{F}^{-T}.
    \end{split}
\end{equation}
Here $\boldsymbol{d}$ denotes the deformation, $\boldsymbol{F}=\frac{\partial x}{\partial X}, \ x \in\Omega^s_t, \  X \in \Omega^s_0 $ the deformation gradient, ${J=det\boldsymbol{F}}$ and the second Piola-Kirchhoff tensors is denoted with  $\boldsymbol{S}$. Outer forces acting on the volume are incorporated in $\boldsymbol{f}_s$.
The elasticity tensor $\boldsymbol{C}=\boldsymbol{C}(E,\nu)$ is given with dependency on the Young's modulus $E$ and the Poisson's ratio $\nu$. The elastic strain tensor is given by the Green-Lagrange strain ${\boldsymbol{\epsilon}=\frac1{2}(\boldsymbol{F}^T\boldsymbol{F}-I)=\tfrac{1}{2}\left(\nabla\boldsymbol{d}+(\nabla\boldsymbol{d})^T+\nabla\boldsymbol{d}\left(\nabla\boldsymbol{d}\right)^T\right)}$. Note, that for small deformation gradients it holds 
$ \boldsymbol{S}\approx \boldsymbol{C:\epsilon}$.\\
The surfaces of the domain clipping, i.e. the annular boundaries of the vessel cuts, surrounding the in- and outflow boundaries of the fluid domain,   $\Gamma_{in}^f$ and $\Gamma_{out}^f$, are constraint to have $0$ deformation. In contrast, the outer surface of the vessel, which would be in contact with surrounding tissue, is able to move freely.\\
To maintain continuity of velocities and forces at the fluid-structure boundary layer, we enforce the coupling conditions to balance velocities and normal stresses of the fluid and the solid material,
\begin{equation} \label{eq:couplings}
    \begin{split}
        \frac{\partial \boldsymbol{d}}{\partial t} &= \boldsymbol{\tilde u}= \boldsymbol{w}\quad \text{and}\\
        J\tilde{\boldsymbol{T}}_f\,\vec{n_f}&= -(\boldsymbol{FS})^T \vec{n_s} \quad \text{on }\Gamma^{fsi}_0.
    \end{split}
\end{equation}
Here, $\boldsymbol{\tilde{u}}, \, \tilde{\boldsymbol{T}}_f$ stand for fluid quantities transformed to the reference fluid-solid layer.\\
Note that the model considering rigid walls consists only of the fluid sub-problem (\ref{eq:ns}) defined in the carotid lumen $\Omega^f=\Omega^f_0, $ moreover ${\boldsymbol w}=0$ and $\frac{D{\boldsymbol u} }{dt}=\frac{\partial{\boldsymbol u} }{\partial t}$. In analogy to the velocity continuity  condition in (\ref{eq:couplings}) the no-slip condition  $\boldsymbol{u}=0$  is prescribed on the vessel wall surface denoted by $\Gamma_w^f$ in case of rigid walls.

\section{Hemodynamic indicators}
\subsection{Wall shear stress} 
The wall shear stress (WSS) measures the endothelial stress 
exerted by blood on the vessel tissue.
To explain the relationship between WSS and zones susceptible to atherosclerosis two main explanatory approaches can be found in the literature. The {\it high shear stress theory} identifies sites with prolonged high WSS as risk zones, the {\it low shear stress theory} considers also sites with oscillating and low WSS as potentially at risk,  the correlation was reported in \cite{ku}. For a systematic review of both see  \cite{pfeiffer}. 

Principally, the wall shear stress 
is  defined on the interface boundary $\Gamma_t^{fsi}$, or on the rigid vessel wall $\Gamma_w^f$, as the projection of the normal stress vector $\vec{\bf t}_f= -{\bf T}_f \vec{n}_f$ onto the tangential plane,
\begin{equation}\label {wss0}
\vec{\tau_w}=\vec{\bf t}_f-(\vec{\bf t}_f\cdot \vec{n}_f)\vec{n}_f=(\vec{\bf t}_f \cdot \vec{t_1})\vec{t_1}+(\vec{\bf  t}_f \cdot \vec{t_2})\vec{t_2},
\end{equation} 
where $\vec{t_1},\,\vec{t_2} $ are unit vectors spanning the tangential plane.  
Alternatively a non-directional quantity describing the amplitude of the wall shear stress vector is frequently evaluated \cite{quarteroni04} as
\begin{equation}\label{wss1}
\tau_w^a=\|\vec{\tau_w}\|=\sqrt{(\vec{\bf t}_f \cdot \vec{t_1})^2+(\vec{\bf  t}_f \cdot \vec{t_2})^2}\,.
\end{equation} 
\\ Different direction-based indicators of WSS have been used  to measure the stress exerted by the fluid  as well, \cite{antiga, gallo, hoogendoorn, mohamied16, mohamied14}. 
For cylinder-like or other simple geometrical objects vector quantities such as the rotary {(transversal)} or longitudinal component of the WSS can be considered.
In  this study, we evaluate the longitudinal component of the wall shear stress vector (longitudinal WSS), aiming to track the backward flow 
in the carotid artery bifurcation vessel tree, which is defined as 
\begin{equation}\label{wss2}
\tau_w^\ell=\vec{\bf t}_f \cdot \vec{t}_\ell=-{\bf T}_f \vec{n}_f\cdot \vec{t}_\ell. 
\end{equation} 
Here, $\vec{t}_\ell$ is a vector of the tangential plane pointing in the longitudinal direction, i.e., the main flow, called "longitudinal tangent" in what follows.    
In (bi-)directional evaluations of WSS, i.e. considering either its longitudinal or transversal component, the proper choice of tangential vector fields on complex surfaces
 is crucial for its evaluation, as it is the case in (\ref{wss2}) for longitudinal WSS.  

\begin{figure}[ht]
\begin{center}
\begin{overpic}[scale=.35,,tics=10]
{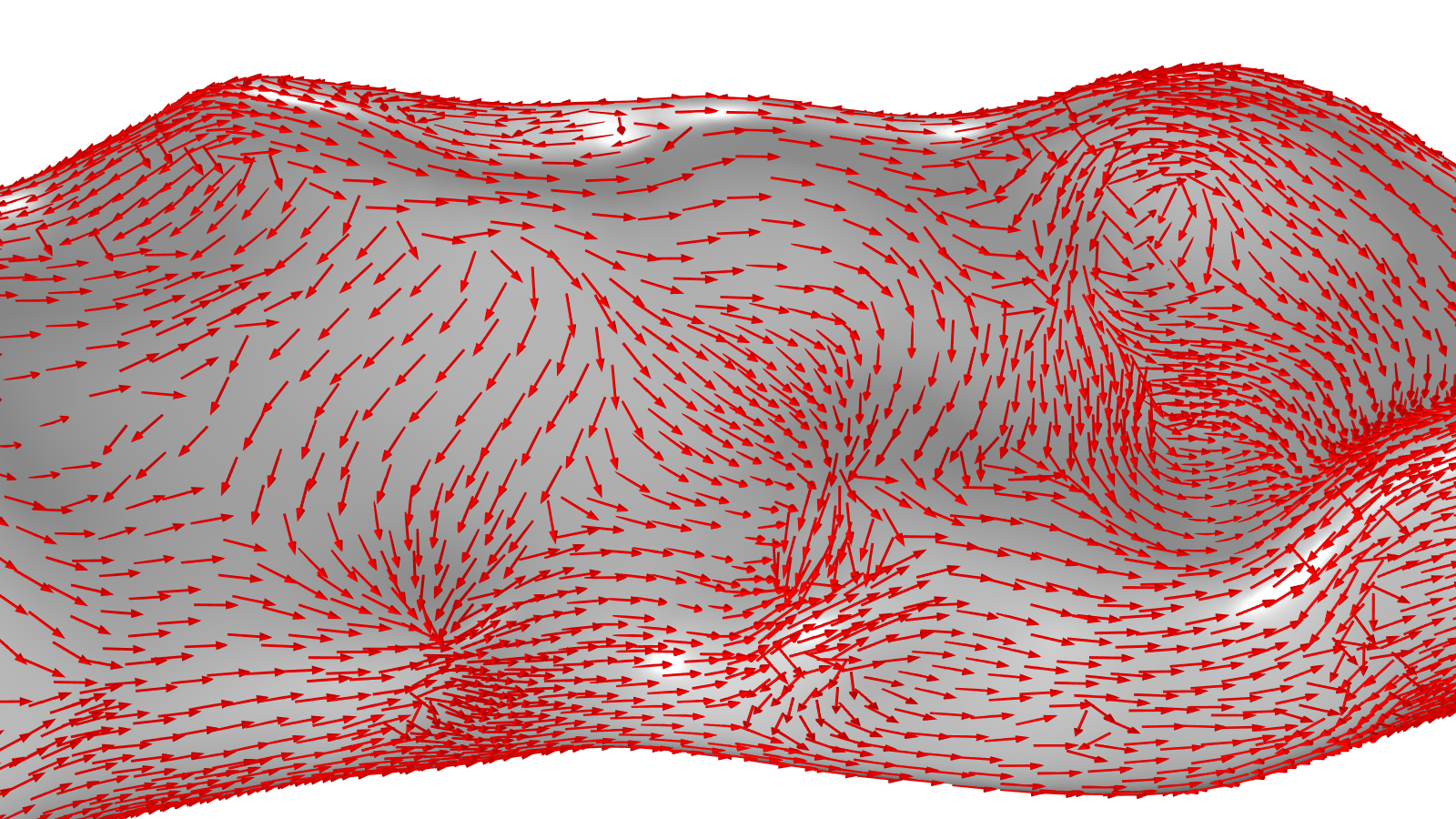}\ \vspace{-2mm} 
\put(10,20){\linethickness{0.25mm}\color{black}\polygon(0,0)(25,0)(25,20)(0,20)}
\end{overpic}
\begin{overpic}[scale=.35,,tics=10]
 {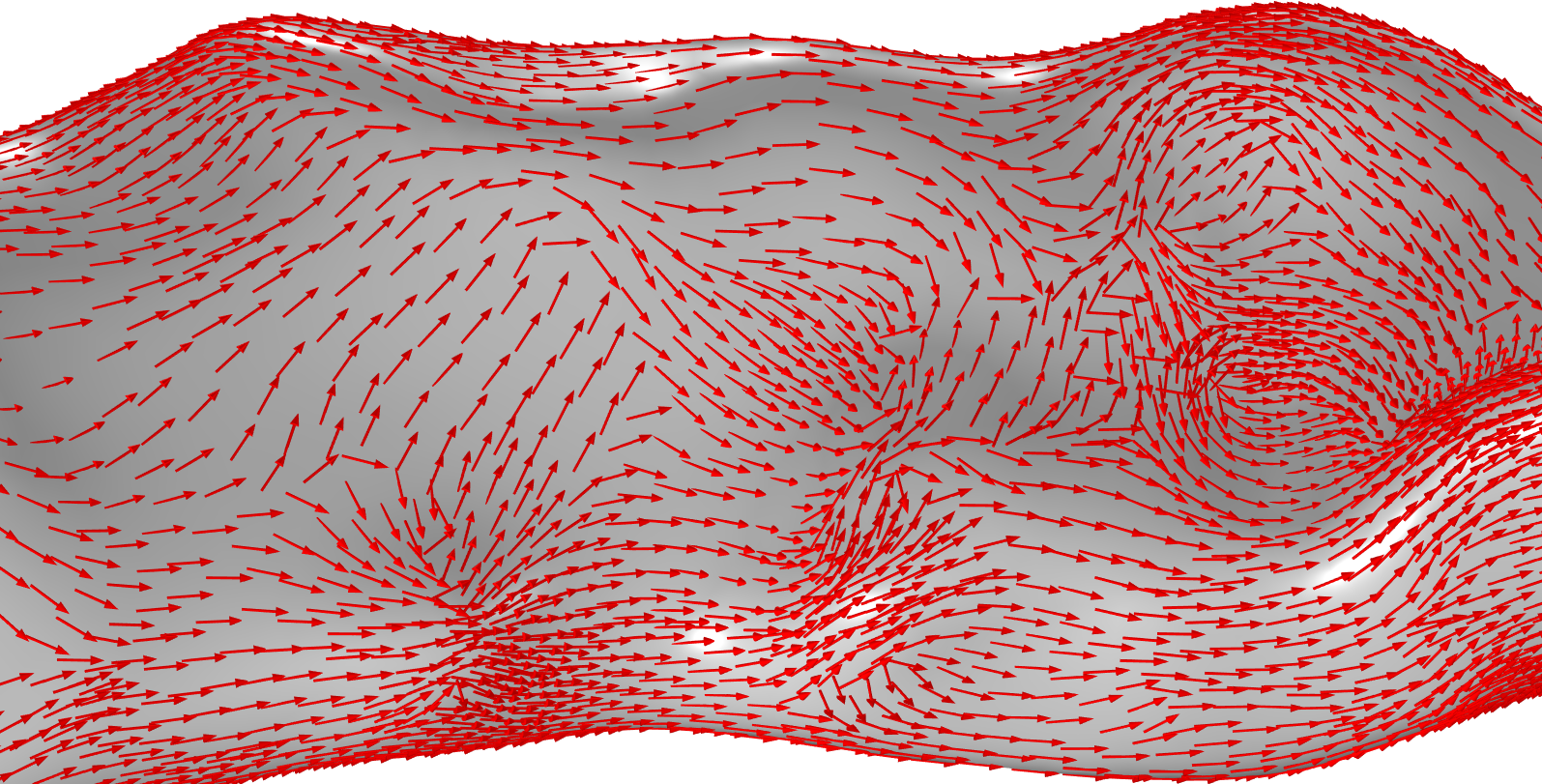}
\put(10,20){\linethickness{0.25mm}\color{black}\polygon(0,0)(25,0)(25,20)(0,20)}
\end{overpic}
\end{center}
\captionsetup{margin=0.5cm}
\caption {Tangential vector field  $\vec{t_2}$ obtained automatically in {\sl Comsol Multiphysics} (left) and flipped tangential field  (\ref{flipped_tang}) (right) on the carotid surface, flow direction: from left to right.}\label{FIG:tangent_unhandled}
\end{figure}

Due to the non-smooth surface of the studied 3D computational geometry, the specification of the proper longitudinal tangent vectors may be problematic. The direction of both tangent vectors $\vec{t_1},\vec{t_2}$ spanning the tangential plane, obtained automatically in {\sl Comsol} for the surface topology,
on neighboring mesh element surfaces may jump and could very well point to the opposite direction of the main flow. This behavior can be observed for $\vec{t_2}$ on the left in Fig. \ref{FIG:tangent_unhandled}. 

In our first approach the tangent vector $\vec{t}_2$ shown in Fig. \ref{FIG:tangent_unhandled} has been chosen as the longitudinal tangent vector $\vec{t}_\ell$, since it fits the longitudinal direction better than $\vec{t}_1$.  However, its discontinuous spatial behavior and its opposing direction to the main flow at some areas would affect the value of longitudinal WSS (\ref{wss2}) substantially. In order to correct the orientation of $\vec{t_2}$ we change it to its diametrically opposite vector. For that we switched its sign according to the angle $\theta$ between $\vec{t_2}$ and an overall flow vector $\vec{v}$. For $\theta \in (\pi/2, 3\pi /2)$ it holds ${\vec{t_2}\cdot \vec{v}} \propto \cos(\theta)<0$. Thus, we define the tangential vector to be used in (\ref {wss2}) as 
\begin{equation}\label {flipped_tang}
    \vec{t}_\ell=\vec{t_2}\,\sign ({\vec{t_2}\cdot \vec{v}}).
\end{equation}
The overall flow  direction vector $\vec{v}$  has to be specified locally for different sections of the computational domain tree. 
 
Note, that $\vec{t}_\ell$ defined by (\ref{flipped_tang}) is revolved from $\vec{v}$ by less then $\pi/2$ and thus aligned almost with the main flow, but it still has the same jumping behavior as the vector $\vec{t}_2$, compare Fig. \ref{FIG:tangent_unhandled} (right). In what follows we present an improved approach for constructing a proper longitudinal tangential field $\vec{t}_\ell$, which is based on the alignment of the carotid tree centerline and can therefore be utilized globally. 

\subsection{Projection method for tangential field}
  
\begin{figure}[H]
  \centering
  \begin{minipage}[b]{.45\linewidth} 
      \includegraphics[width=\linewidth]{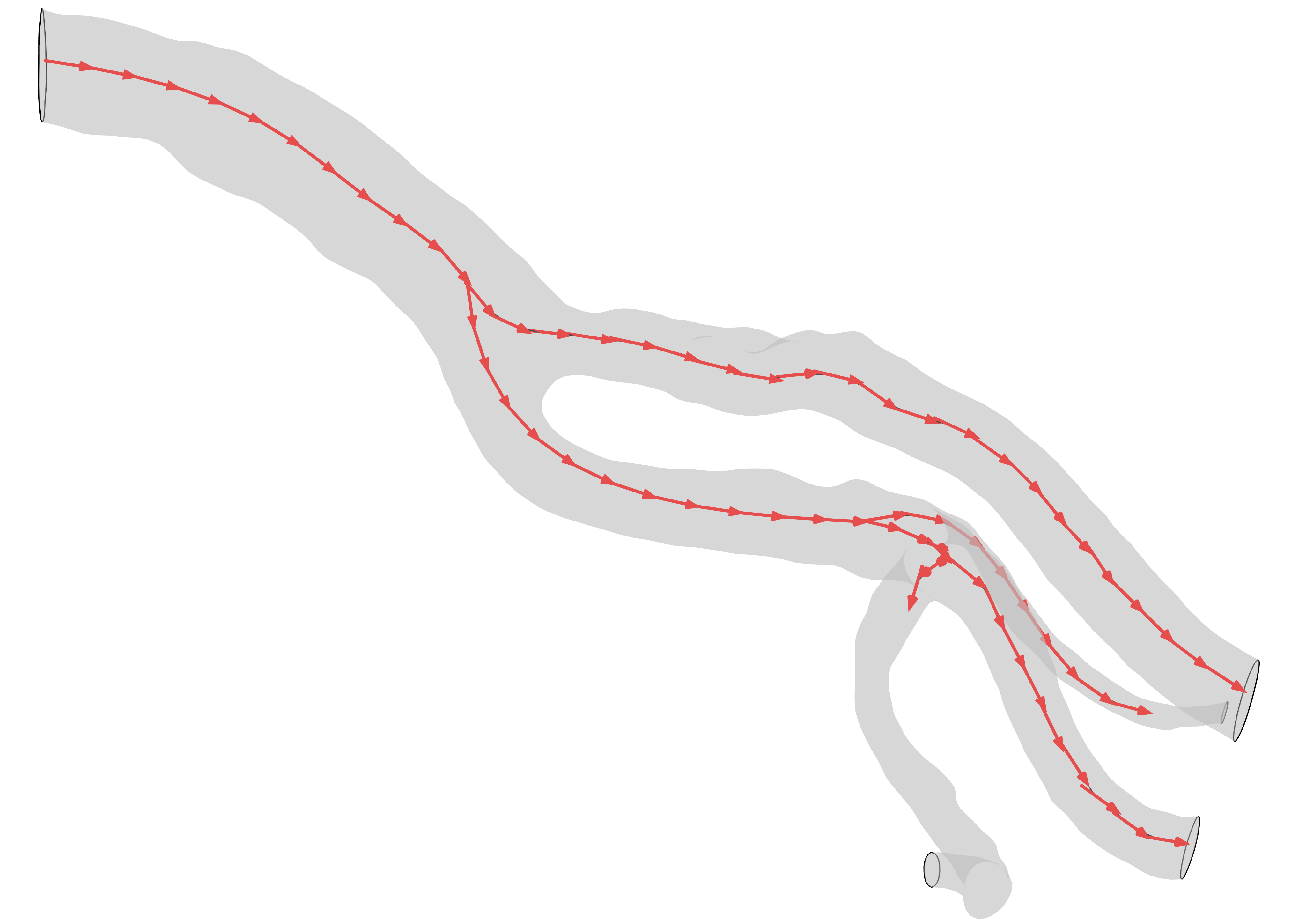}
      \caption{The carotid domain and the normed tangential field  $\vec{c}^{\,\prime} $ of the centerline.}
     \vspace{-0.5cm}
     \label{Fig:centerline tangential fiel}
  \end{minipage}
  \hspace{.05\linewidth}
  \begin{minipage}[b]{.45\linewidth} 
      \includegraphics[width=\linewidth]{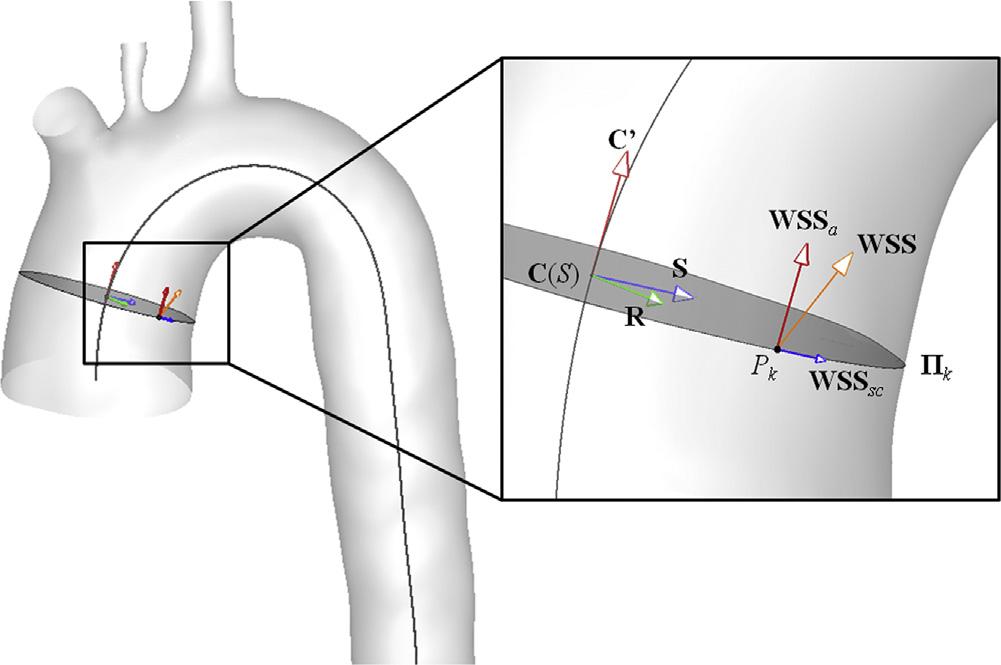}
      \caption{Projection method by Morbiducci, adapted from \cite{morbiducci15}.}
      \vspace{0.45cm}
      \label{Fig:projection method}
  \end{minipage}
\end{figure}

 As depicted above in Fig. \ref{FIG:tangent_unhandled}, on complex surfaces the  automatically rendered  tangent vectors $\vec{t_2}$ do not follow the overall flow direction in some topologically complicated areas. In order to overcome this difficulty we apply  the approach based on the knowledge of the centerline  of the vessel tree and its projection onto the  vessel surface, similarly to the method of Morbiducci et al \cite{morbiducci15}, presented in Fig.~\ref{Fig:projection method}.


At first the centerline is obtained as the set of center-points of the maximally inscribed spheres.
Here, we use the 3D Voronoi diagram of the geometry to find and connect the center-points, the method described and  implemented in the vascular modeling toolkit~\cite{izzo2018}.
The method yields robust and detailed results
with a resolution of about 3000 points.
 After getting the 3D curves of the centerline, its tangential vectors $\vec{c}$ are obtained by subtraction of  two points of the curve, see Fig.~\ref{Fig:centerline tangential fiel}.

\noindent

\vspace{0.3cm}

\begin{figure}[H]
\centering
\begin{minipage}[b]{.60\linewidth}
    \includegraphics[width=\linewidth]{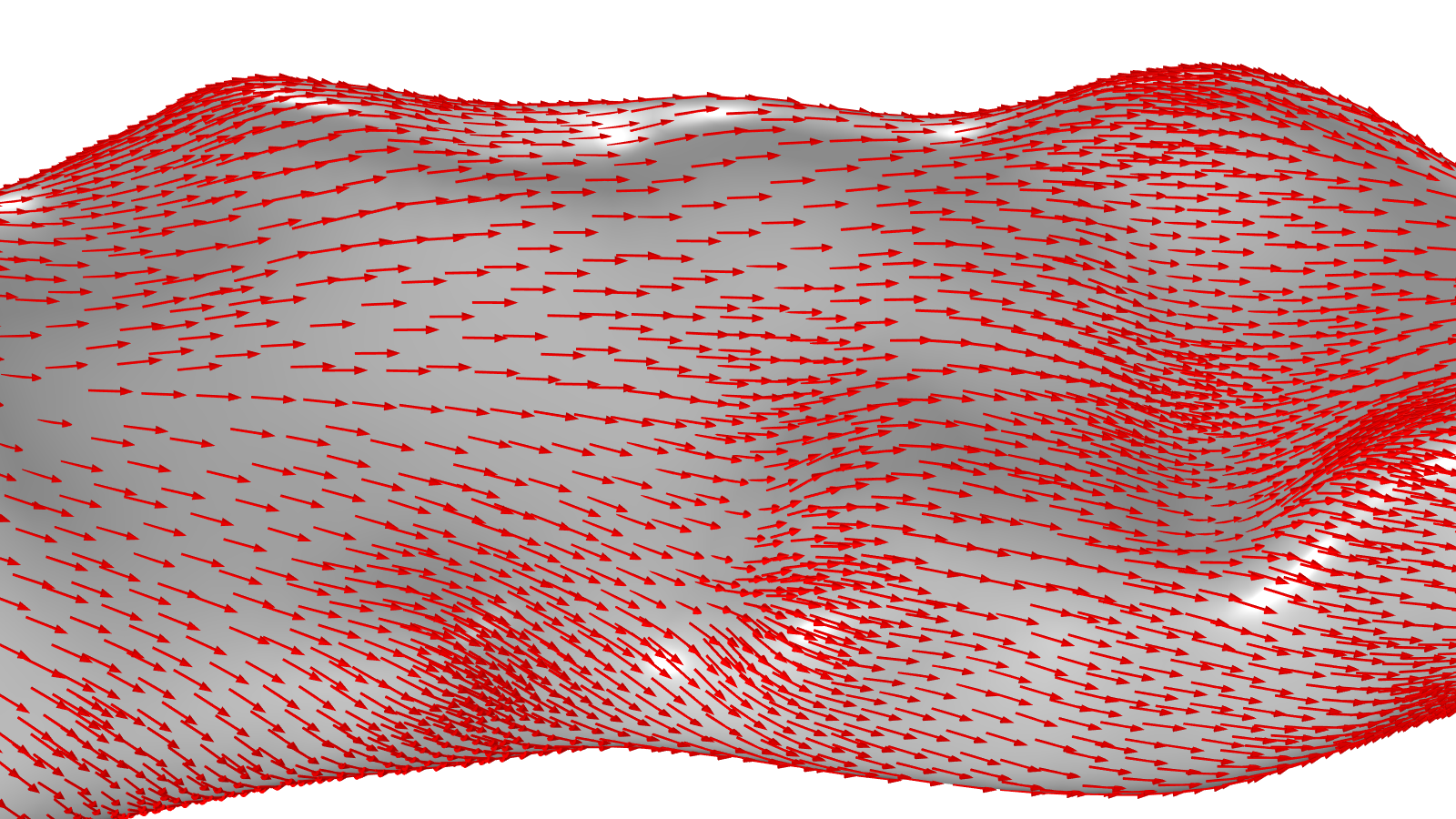}
    \caption{Tangential vector field $\vec{t}_\ell$ obtained from the centerline by the projection method.}
    \label{fig:project_tang} 
\end{minipage}
\end{figure}

 Afterwards, the centerline tangent vectors $\vec{c}^{\,\prime} $ are projected to the vessel surface, i.e., into each surface point $P_k$. This is done in two steps, first $\vec{c}^{\,\prime}$ is extrapolated to the surface points by the geometry tool {\sl extrapolate} with {\sl linear} settings in {\sl Comsol}. Then, the extrapolated centerline tangents $\vec{c}$ are projected into the tangential plane of the carotid artery surface by subtracting its normal component, 
\begin{equation}\label {proj_tang}
\vec{t}_\ell= \frac{\vec{c}-(\vec{c}\cdot \vec{n}_f)\vec{n}_f }{\| \vec{c}-(\vec{c}\cdot \vec{n}_f)\vec{n}_f\|},
\end{equation} 
here $\vec{n}_f$ are the normal vectors of the carotid surface. The procedure is illustrated in Fig. \ref{Fig:projection method}, where the centerline tangent is denoted by $C'$ and the corresponding longitudinal component of the WSS by WSS$_a$.
The resulting longitudinal tangential field $\vec{t}_\ell$, see Fig. \ref{fig:project_tang}, shows a more uniform alignment at first sight 
compared to the flipped tangential field (\ref{flipped_tang}) presented in Fig. \ref{FIG:tangent_unhandled}.
In this manner, longitudinal tangent vectors $\vec{t}_\ell$ derived from the centerline, with proper unidirectional behavior on the surface, are implemented in {\sl Comsol}. Both, projected (\ref{proj_tang}) as well as flipped tangents (\ref{flipped_tang}) are used for the evaluation of longitudinal WSS and the results are compared in what follows.

\subsection{Oscillatory shear index}\label{subsec:osi}
The oscillatory shear index (OSI) introduced by Ku et al. \cite{ku} is a common indicator for disturbed flow. It characterises
the temporal oscillations of  WSS through its directional change at any point on the surface in the considered time period. The  degree of oscillation is expressed by the 
ratio  of averaged WSS compared to its averaged  amplitude over the whole time interval, i.e in means of temporal mean values. Note, that OSI does not express the frequency of the sign change of the WSS.
\\ We introduce two definitions of OSI, which  can be found in literature, see e.g., \cite{antiga, blagojevic, Qua_OSI, quarteroni04, soulis, Tremmel_OSI, taylor}, based either  on the amplitude of the temporal mean of the WSS vector, (\ref{wss1}), or on the size and sign of the temporal mean of its longitudinal component (\ref{wss2}),
\begin{equation}\label {osi2}
OSI=\frac1{2}\left ( 1- \frac{\|\int_0^T \vec{\tau_w}\,\mathrm{dt}\,\|}{\int_0^T \|\vec{ \tau_w}\| \mathrm{dt}}\right),  
\end{equation} 
\begin{equation}\label {osi1}
OSI^\ell=\frac1{2}\left ( 1-\frac{\int_0^T \tau_w^\ell \  \mathrm{dt}}{\int_0^T |\tau_w^\ell| \mathrm{dt}}\right).  
\end{equation} 
 Note, that these formulas differ in the ratio and the norm of the ratio of temporal mean WSS. Consequently, formula (\ref{osi2}) defines values of the OSI between 0 and 0.5, where 0 stands for no or a complete change of sign over the mean WSS's entire time interval 
 and 0.5 for completely balanced sign changes, e.g., oscillations  of the mean WSS. Values in between 
 imply corresponding sign changes of the WSS.
 In contrast, negative values of the mean longitudinal WSS in (\ref{osi1}) play an important role. They lead to a range of OSI values from 0 to 1. Similarly to definition (\ref{osi2}), values of 0.5 indicate completely balanced sign changes of WSS as well. An OSI value of 0 represents a point at which the WSS is positive over the entire time interval considered. For the value 1, on the other hand, the WSS is negative over the entire time interval. Values between 0 and 0.5 show predominantly positive, values between 
 0.5 and 1 describe predominantly negative  WSS over the whole time interval.
 The definition (\ref{osi1}) thus allows to locate not only sites of oscillating WSS, but also sites with long-lasting or predominantly negative WSS. Thus, in contrast to  (\ref{osi2}),  (\ref{osi1}) provides an index that can represent both indicators of \textit{low shear stress theory}. In what follows, 
 we refer to the directional definition (\ref{osi1}) when mentioning the OSI, 
 but we also evaluate the OSI defined with the use of the WSS amplitude (\ref{osi2}).
 
 \section{Numerical method and convergence study}
 The numerical simulations have been performed with {\sl Comsol Multiphysics} \cite{comsol}, Version 5.6 using the {\sl MEMS Module} to incorporate the interaction of laminar fluid flow with the linear elastic wall material. The software uses the arbitrary Eulerian-Lagrangian formulation, which consists of the Navier-Stokes equations in the Eulerian frame (\ref{eq:ns}) and the solid mechanics equations using the Lagrangian formulation (\ref{eq:def}). The interaction was chosen to be bidirectional, such that the fluid loading acted on the structure and the wall velocity is transmissioned to the fluid. We chose a monolithic, fully-coupled solving scheme to provide a robust solution for the dependent variables consisting of the solid deformation ${\boldsymbol d}$, the fluid flow velocity ${\boldsymbol u}$, pressure $p$ and the spatial mesh displacement ${\boldsymbol w}$.  
 
\begin{figure}[H]
    \centering
    \includegraphics[scale=0.5]{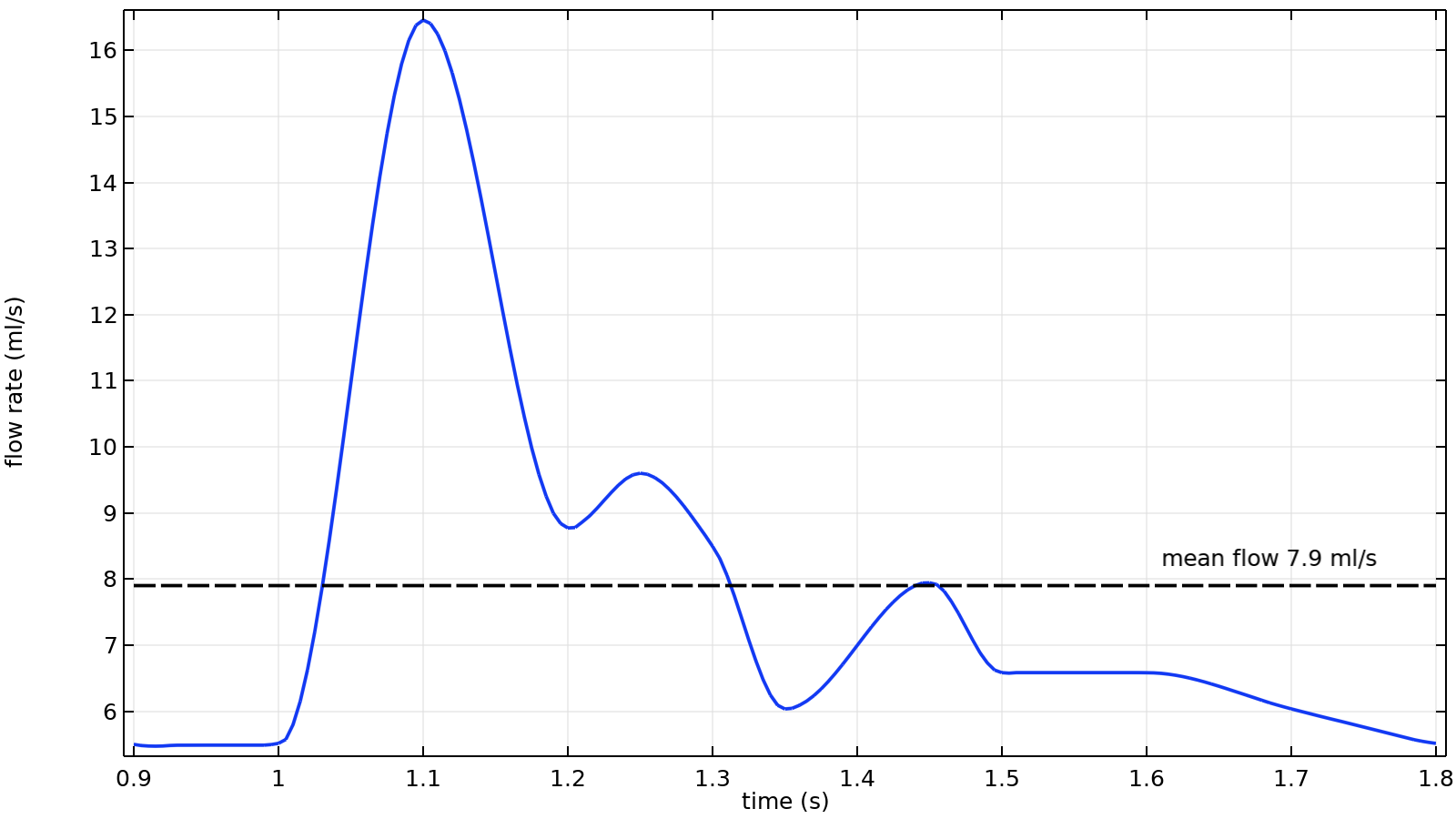}
    \caption{Flow waveform at the inflow boundary of the common carotid artery with 7.9 ml/s mean flow rate, waveform adapted from \cite{antiga, perktold95}.}
    \label{fig:flow rate}
\end{figure}

 For the discretization of the fluid  velocity and pressure linear 
and for the solid deformation  quadratic 
finite elements have been chosen. For the time discretization BDF-method 
of order 1 (for the initialization) continuing with order 2 and an adaptive time-stepping has been applied. The simulations have been performed on a computational mesh consisting of about 
$163500$
tetrahedral elements in the solid domain and 
$510800$
elements in the fluid domain,  about 
$412500$
of which are tetrahedral and 
$98300$ are prisms acting as two boundary layers.
The whole simulation spans over two cardiac cycles, i.e. $1.8s$, and is driven by a pulsatile flow rate presented in Fig. \ref{fig:flow rate}, starting with the domain at rest, i.e.,  zero deformation and zero flow at $t=0$. All  results presented in  the next section are chosen from the second cycle where the flow is fully developed. 

\subsection {Convergence study}
\begin{figure}[h!]
  \begin{minipage}[b]{.49\linewidth} 
       \includegraphics[width=\linewidth]{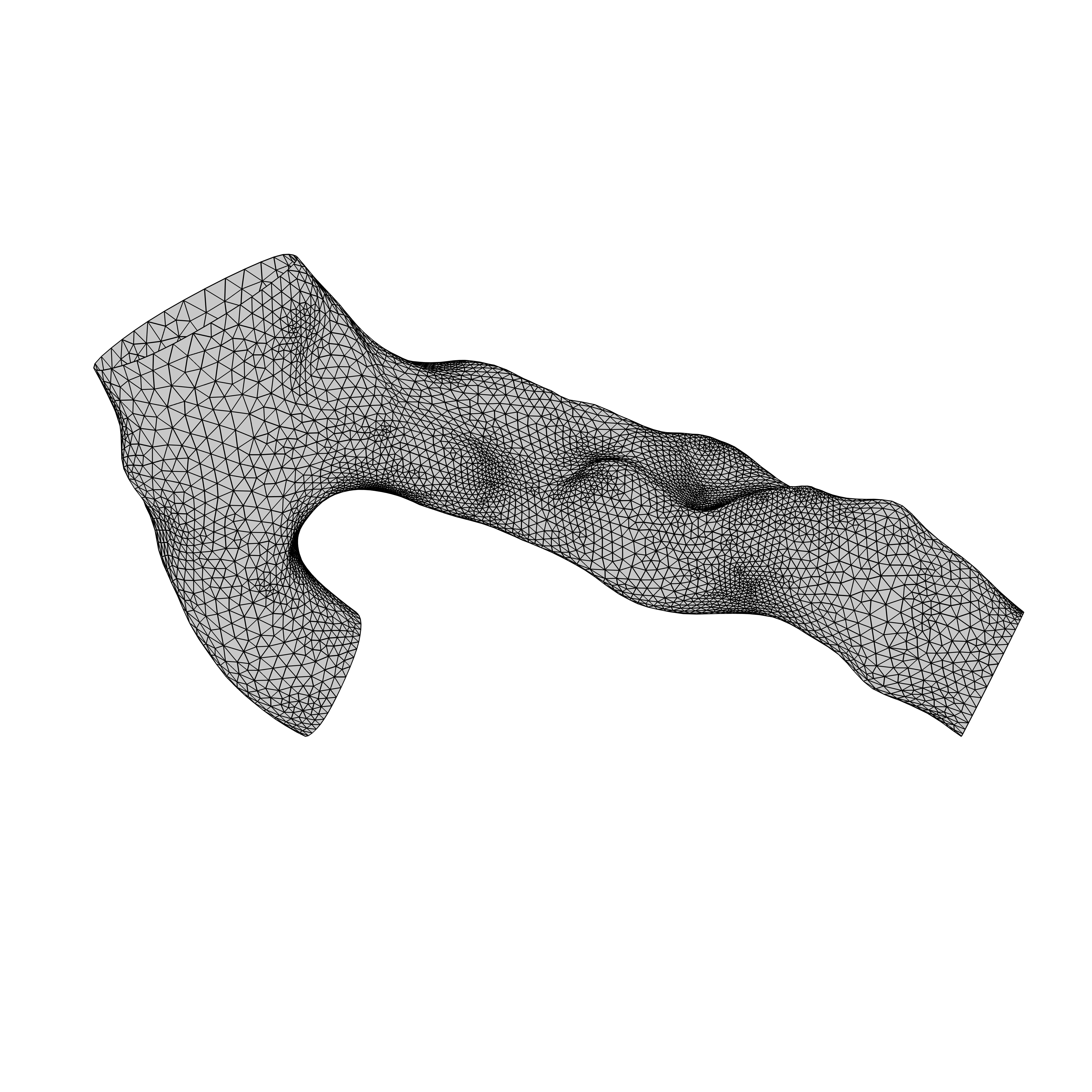}
      \caption*{Solid wall model mesh}
  \end{minipage}
  \hspace{.05\linewidth}
  \begin{minipage}[b]{.49\linewidth} 
      \includegraphics[width=\linewidth]{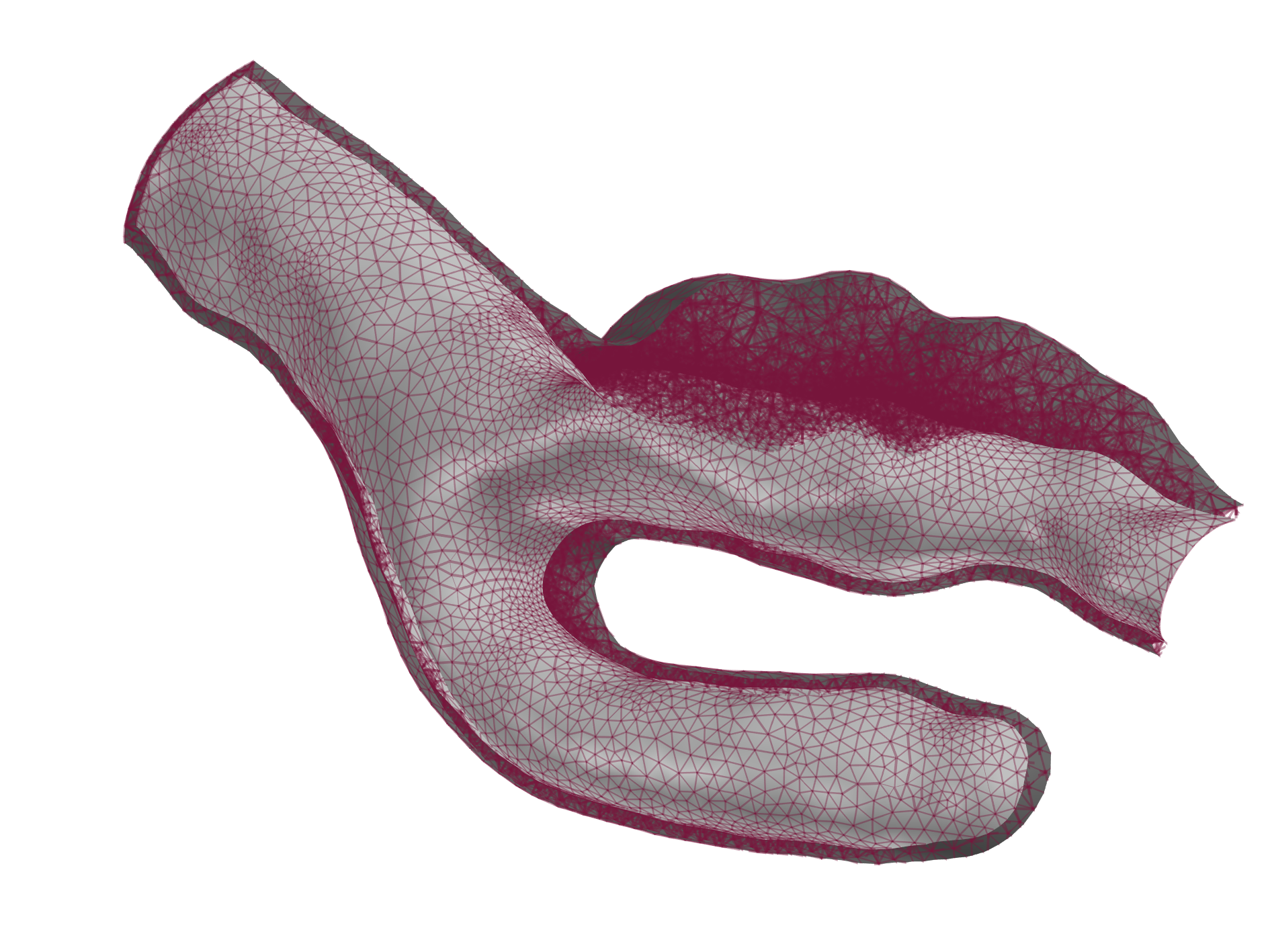}
      \caption*{Mesh for compliant walls model (FSI)}
  \end{minipage}
  \caption{A part of computational mesh for solid and compliant vessel wall model,  chosen sections of interest for data evaluation presented.}
      \label{Fig:mesh}
\end{figure}

The numerical mesh convergence study has been performed for the fluid flow problem with solid vessel walls to restrict  the computational costs. The spatial mesh error was computed using a set of eight meshes, approximately doubling the mesh element number from  mesh $i$ to mesh $i+1$. To compare numerical solutions on different meshes
with non-coinciding mesh nodes, the linear shape functions (P1 finite elements) are applied. The solution difference is realised with use of {\sl join solution} feature in {\sl Comsol} by projection of the lower-mesh solution onto the higher mesh.  The spatial discretization error has been evaluated  for fluid velocities and for the longitudinal WSS in means of weighted $L^2$- norm of the difference of the reference and actual $i$-th mesh solution obtained on meshes Nr. 1-6, 
\begin{equation}
    \begin{split}\label{eq:errors}
    err(\boldsymbol{u} _i)&:=\frac1{\sqrt{|\Omega^f|}}\|\boldsymbol{u} _i-\boldsymbol{u}_{ref}\|_{L^2(\Omega^f)}, \\
    err(\tau_i)&:=\frac1{\sqrt{|\Gamma_w^f|}}\|(\tau_w^\ell)_i - (\tau_w^\ell)_{ref}\|_{L^2(\Gamma_w^f)}, \quad i=1, \ldots 6,
    \end{split}
\end{equation}
in the time-point of the maximal flow of the second cardiac cycle, $t=1.1$ s.
Here, $\boldsymbol{u}_{i},\ (\tau_w^\ell)_i$ are numerical data obtained on the $i$-th mesh. The reference  solution has been obtained on mesh Nr.  8 
consisting of a total of 8 717 584 tetrahedra and prism elements. The mesh errors (\ref{eq:errors}) have been evaluated in the chosen section of our computational mesh presented in Fig. \ref{Fig:mesh} (left), which is identical to the geometrical region of all numerical results presented.
 
\begin{figure}[H]
    \centering
    \includegraphics[scale=0.27]{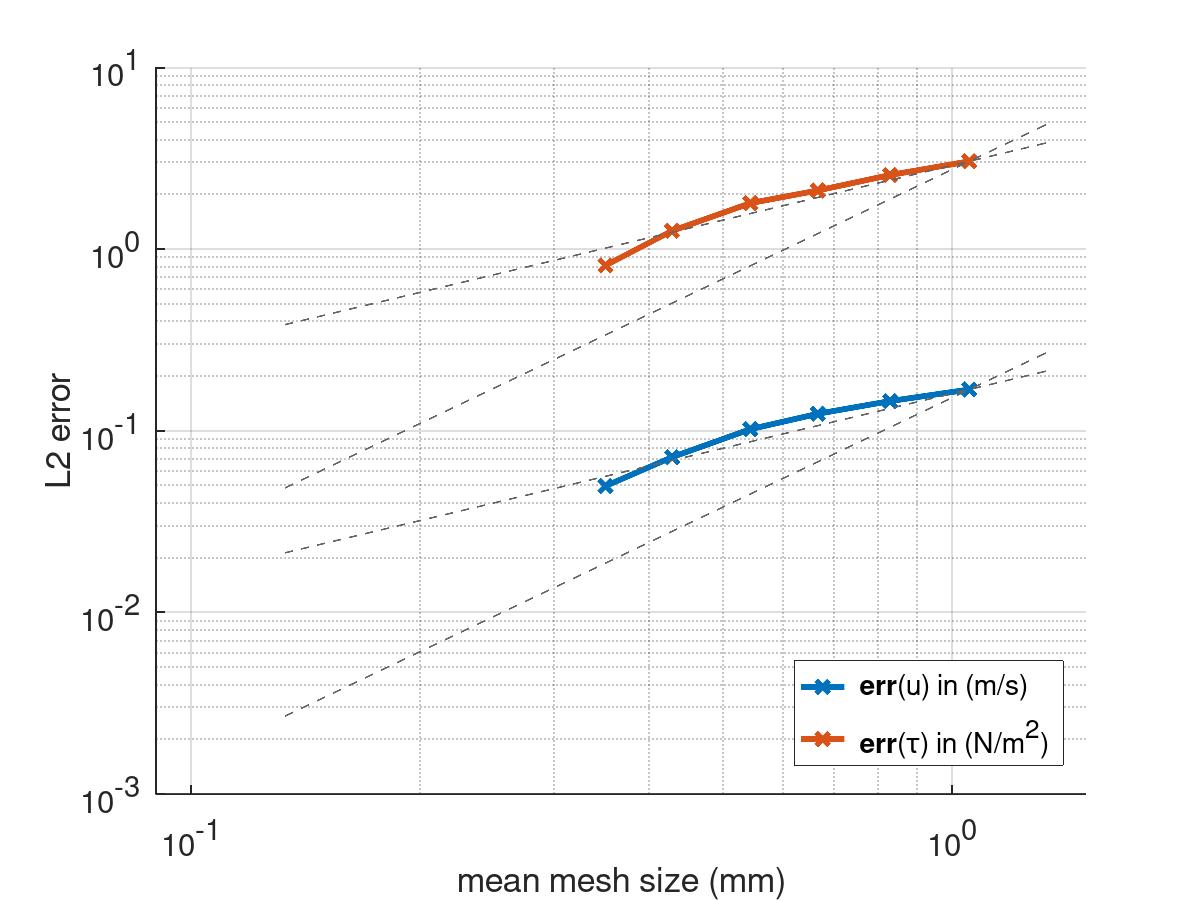}
    \caption{Decrease of $L^2$ error norm (\ref{eq:errors}) for $\boldsymbol{u}$ and  $\tau_w^\ell$ with respect to the mean mesh size, gray dashed lines with slopes 1 and 2 represent first and second converegence order.}
    \label{fig:convergence}
\end{figure}

The mesh element sizes, measured by the diameter the maximum inscribed sphere, have been averaged for each mesh and  spatial errors have been related to the mean mesh sizes $h_i$. The descending sequence of $h_i$ starting with  $h_1=1.054$ mm is presented in Table \ref{tab:convergence}, whereby the decrease factor $$a_i:=h_i/h_{i+1}$$  lies  between  1.271 and 1.22.
In Fig. \ref{fig:convergence} the  errors with respect to $h_i$ are presented in logarithmic scale. For the decreasing error curves one can observe a  slope of approximately one, which slightly differs for higher and lower mesh sizes. In analogy to the error estimation for finite element method with a uniform mesh and constant element sizes $h$, $\|u_h-u_{exact}\|_{L^2} \approx Ch^p,\,  C < \infty$,
where the convergence order $p$ is identified with the slope of the logarithmic error curves for different $h$,  we denote the slope  of our logarithmic error curve by the experimental order of convergence (EOC).  EOC can be obtained  by comparing errors of two consecutive meshes, 
\begin{equation} \label{eoc}
EOC (\boldsymbol{u}_i)=\frac{\log_{10} (err(\boldsymbol{u}_i))-\log_{10}(err(\boldsymbol{u}_{i+1}))  }{\log_{10}h_i-log_{10} \, h_{i+1}}
=\log_{a_i}\left( \frac{err(\boldsymbol{u}_i)}{err(\boldsymbol{u}_{i+1})}\right),   
 \end{equation}
analogously for $(\tau_w^\ell)_i, \ i=1, \ldots 6$.

 
 The results of our mesh convergence study are summarised in Table \ref{tab:convergence}  and show good convergence of the  numerical velocities as well as for longitudinal wall shear stresses to  the reference numerical solution. The EOC,  initially lower than one, increases continuously with decreasing mesh size, with a small deviation in case of $(\tau_w^\ell)_3$. The averaged EOC for the velocities is about 1.13 and for the longitudinal WSS it is about 1.21, altogether a  super-linear averaged convergence  to the reference solution is obtained in our simulations  with rigid vessel walls. 
 Applying
  the finite element method with linear P1 elements, second order convergence rate is expected for spatial error. The decreased experimental order of convergence is caused by the discontinuity of the boundary conditions for velocity  on the edges common  for the inflow  $\Gamma_{in}^f$  as well as vessel wall boundary $\Gamma^f_w$, where  the no-slip boundary condition  meets the non-zero  inflow velocity, considered to be constant  over $\Gamma_{in}^f$.
\begin{table}[h!]
    \centering
      \begin{tabular}{|c|c|c|c|c|c|c|}
      \hline
    mesh  $(i)$  & \# of elements & $h^i$ (mm) & $err(\boldsymbol{u}_i)$ (m/s) & $EOC(\boldsymbol{u}_i)$  & $err({\tau}_i)$ (N/m$^2$) & $EOC({\tau_w^\ell})_i$\\ \hline \hline
     1.   & 6.106E+04 & 1.054 & 1.689E-01 & 0.613   & 3.046E+00 & 0.726   \\ \hline 
     2.   & 1.250E+05 & 0.829 & 1.458E-01 & 0.736   & 2.559E+00 & 0.896   \\ \hline 
     3.   & 2.484E+05 & 0.667 & 1.241E-01 & 0.956   & 2.104E+00 & 0.788   \\ \hline 
     4.   & 4.584E+05 & 0.544 & 1.021E-01 & 1.497  & 1.792E+00 & 1.480   \\ \hline
     5.   & 1.049E+06 & 0.429 & 7.161E-02 & 1.819   & 1.262E+00 & 2.179   \\ \hline 
     6.   & 2.079E+06 & 0.351 & 4.965E-02 & - & 8.135E-01 & -   \\ \hline 
     \end{tabular}
    \caption{Spatial discretization errors related to mean mesh size.}
    \label{tab:convergence}
\end{table} 

\section{Results and discussion}
\noindent

\noindent

The evaluations of  numerical data and  considered wall parameters 
are compared for four model configurations. In two different vessel wall models: rigid $\&$ elastic vessel walls two different tangential fields $\vec{t}_\ell$: flipped  
(\ref {flipped_tang}) $\&$ projected tangential field (\ref {proj_tang}) have been used to calculate the longitudinal WSS (\ref{wss2}) and the corresponding $OSI^\ell$ (\ref{osi1}),  the overview of model configurations can be found in Table \ref{table:configurations}. 
Additionally, the amplitude-based wall parameters (\ref{wss1}), (\ref{osi2}) are evaluated as well  and  compared  to its directional counterparts.
\begin{table}[H]
    \centering
    \begin{tabular}{c|c}
      Configuration  & Properties  \\
      \hline
        (a) &  rigid walls / flipped tangents \\
        \hline
        (b) &  rigid walls  / projected tangents\\
        \hline
        (c) &  elastic walls / flipped tangents \\
        \hline
        (d) &  elastic walls / projected tangents\\
        \hline
    \end{tabular}
    \caption{Configurations of the model evaluations for directional wall parameters $\tau_w^\ell$ and $OSI^\ell$.}
    \label{table:configurations}
\end{table}
In every simulation fluid density $\rho_f=1000 $ kg.m$^{-3}$ and constant viscosity $\mu={0.00345}$ Pa.s was chosen. The wall parameters density $\rho_s={1070}$ kg.m$^{-3}$, Young's modulus $E={0.5}$ {MPa} and Poisson's ratio $\nu=0.17  $ $\nu \in \left[0.17,0.5\right)$\footnote{The choice of Poisson's ratio has negligible impact on the vessel wall displacement in our modeling.} were used in simulations (c) and (d) with fluid-structure interaction.

\medskip

\begin{figure}[h!]
    \centering
    \begin{subfigure}[b]{0.4\textwidth}
        \centering \begin{overpic}[height=0.29\textheight,,tics=10]
        {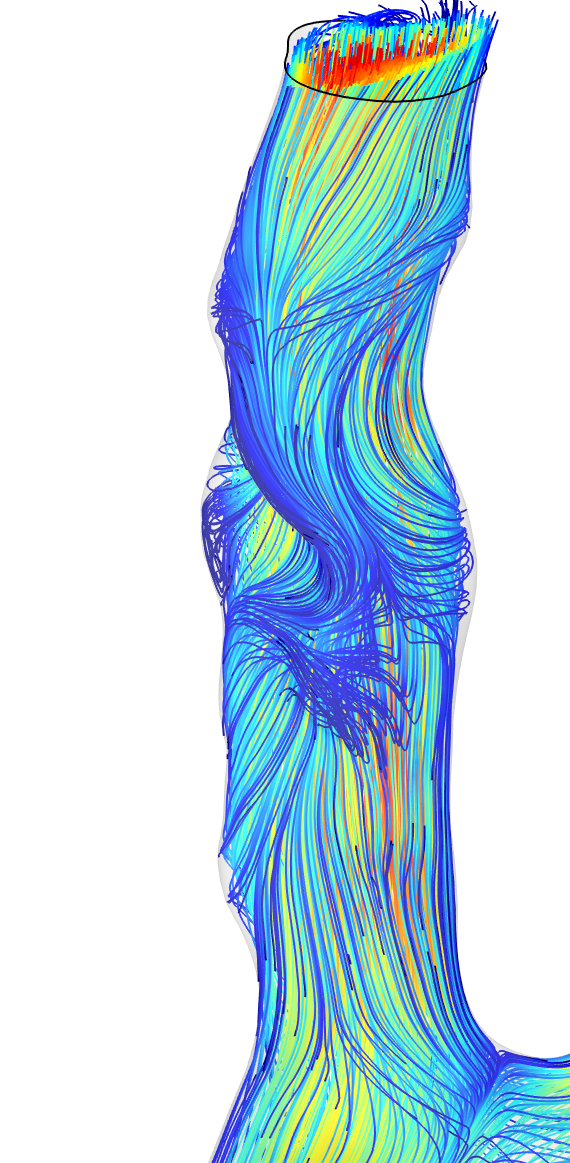}
        \put(9,63){\linethickness{0.25mm}\color{black}
        $\longrightarrow$}
        \put(39,66){\linethickness{0.25mm}\color{black}
        $\longleftarrow$}
        \put(39,63){\tiny stenosis}
        \end{overpic}
        \includegraphics[height=0.29\textheight]{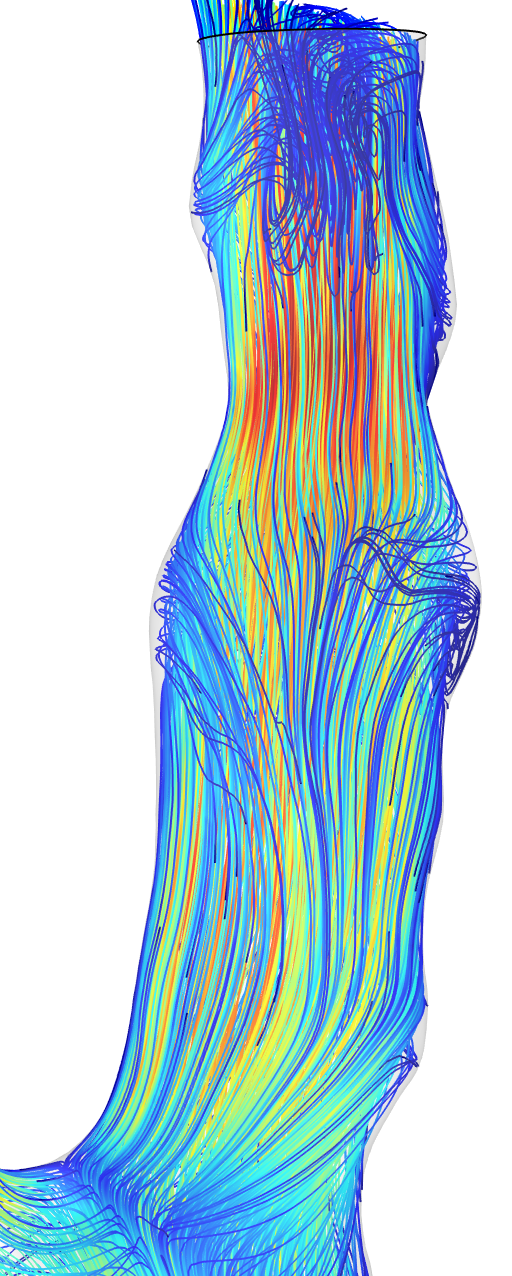}
        \caption*{Rigid walls}
    \end{subfigure}
    \hfill
    \begin{subfigure}[b]{0.18\textwidth}
        \centering
        \includegraphics[height=0.31\textheight,width=0.17\textwidth]{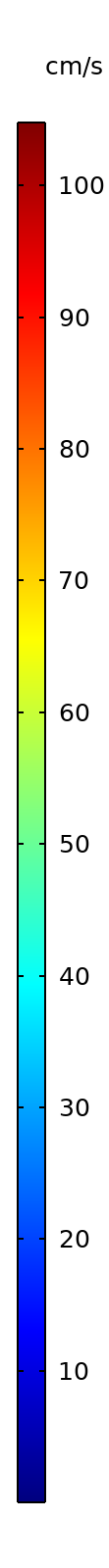}
        \caption*{}
    \end{subfigure}
    \hfill
    \begin{subfigure}[b]{0.4\textwidth}
    \centering
       \includegraphics[height=0.29\textheight]{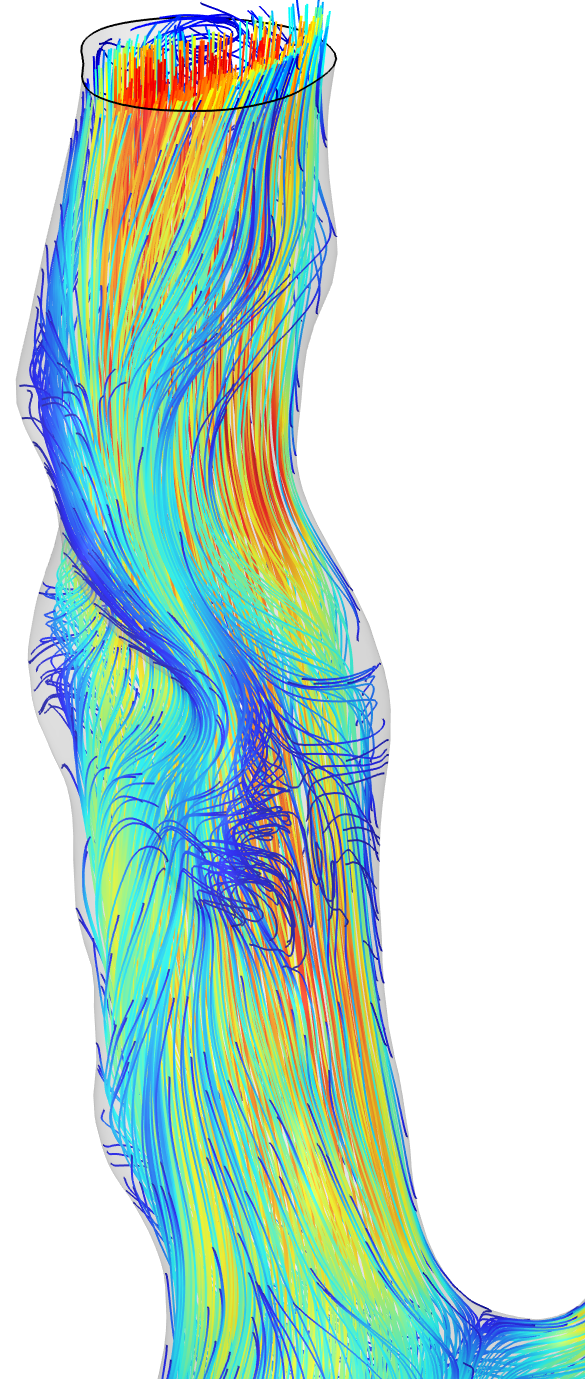}
        \includegraphics[height=0.29\textheight]{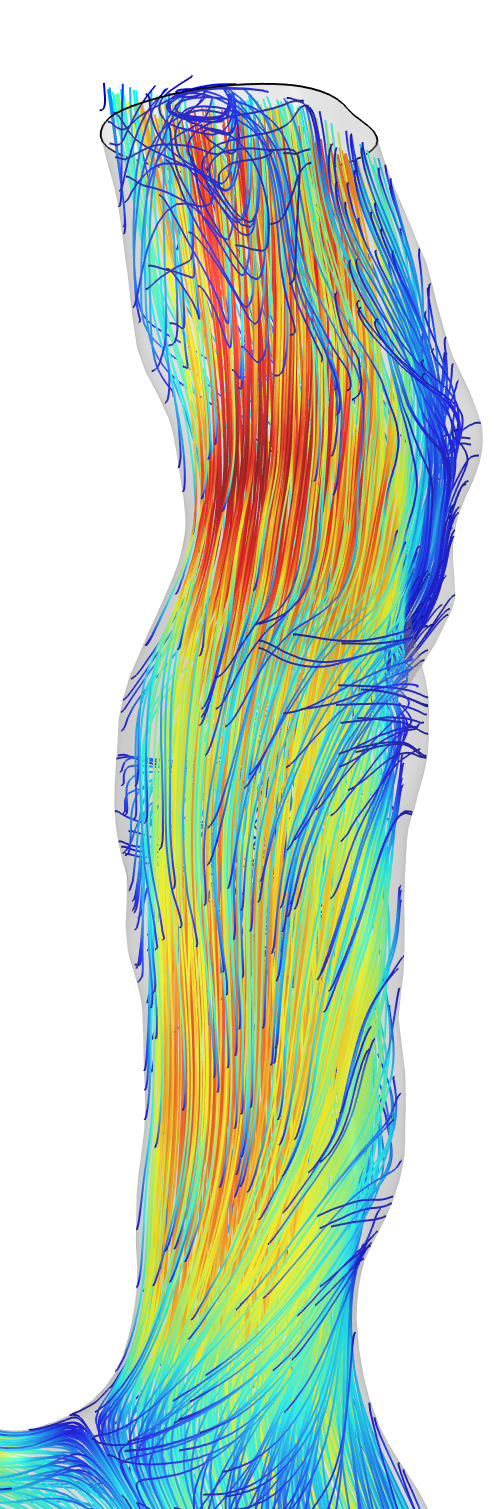}
        \caption*{Elastic walls}
    \end{subfigure}
    \caption{ Velocity streamlines in the stenosed region of ICA (stenotic plaque marked with arrows) for both vessel wall models colored by  velocity magnitude, $t=1.1$ s,  viewpoints from front and back.}
    \label{fig:streamlines}
\end{figure}

In Fig. \ref{fig:streamlines} the blood velocity streamlines in the stenosed region of the ICA are presented {at the time  of maximal flow-rate, $t=1.1s$}. The observed vortices are located around the stenotic bulges in the area adjacent to the stenosis. Note, that the appearance of vortices and backward flow is related to high $OSI^\ell$ values, observed in Fig. \ref{fig:OSI}, and may imply progression of lesions along the carotid artery tree, which is a common hemodynamic hypothesis \cite{spanos}. Note, that amplitude-based OSI presented in Fig. \ref{fig:osi2} only indicates the edges of these areas. On the other hand, high-valued and unidirectional velocity streamlines are observed along the inner wall of the ICA and are related to high longitudinal WSS values $\tau_w^\ell$, as well as the amplitude $\tau_w^a$  of WSS vector, observed  in Figs. \ref{fig:wss_surface}, \ref{fig:WSS_vector}. The latter fits to the velocity profile observations in \cite{zarins_streamlines}.

\begin{figure}[h!]
  \centering
      \includegraphics[width=0.6\linewidth]{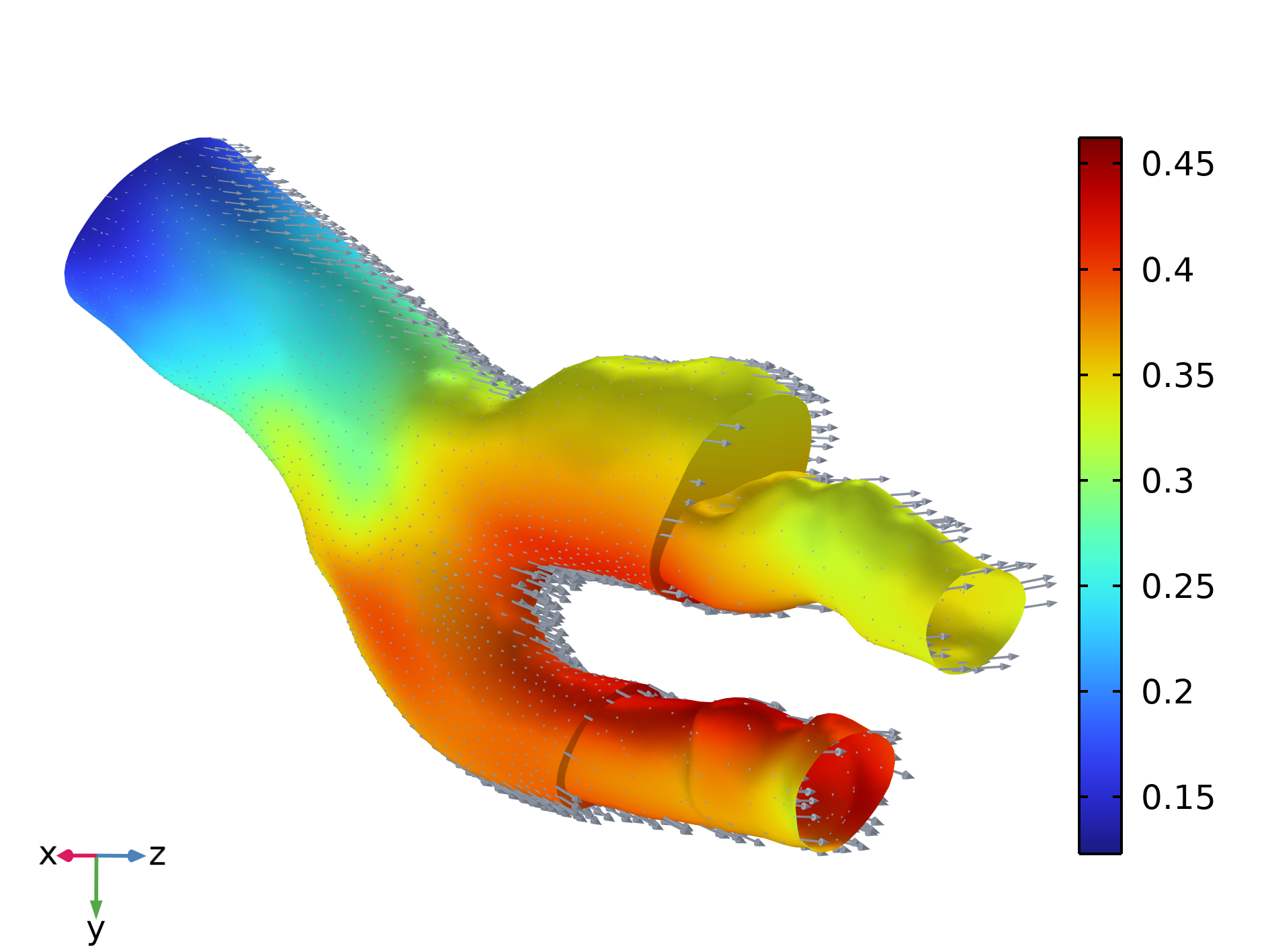}
      \caption{The deformation vectors at the inner and outer surface of the thick vessel tissue shows the translational vessel movement in the zoomed area prior and subsequent to the carotid bifurcation. The surface is coloured by the deformation amplitude in cm.}
      \label{Fig:wall_deformation}
\end{figure}

In Fig. \ref{Fig:wall_deformation} the movement of the carotid in regions close to its bifurcation is presented, whereby the clipped ends of the vessel are cut out. The arrows at the common fluid-solid interface (inner wall) and at the outer vessel wall show the dominance of the lateral vessel wall movement compared to the inflating effects. On the one hand, this is caused by the  boundary condition at the outer wall, allowing free movement without any constraining effects of outer stresses by surrounding tissue, on the other hand, by the length of the whole computational domain considered. Indeed, the effects of clipping the bottom and top part of the vessel geometry, whose length is identical to the fixed geometry presented in Fig. \ref{Fig:Computational Geometry}, are weakened in the middle of the computational domain. This leads to the fluid-driven lateral motion. 

\subsection{Wall shear stress}
In what follows we present the wall shear stress computed from numerical data in means of its longitudinal component, $\tau_w^\ell$, (\ref{wss2}),  as well as the non-directional quantity $\tau_w^a$, (\ref{wss1}) measured by the amplitude of the WSS vector. To compare the WSS distributions for compliant and rigid walls the results are presented on the inner arterial wall in the reference  geometry frame.
\subsubsection{Longitudinal WSS}
\noindent
The longitudinal  WSS (\ref{wss2}) is evaluated on the carotid surface as well as along chosen surface curves and compared for model configurations (a)-(d), see Table \ref{table:configurations}.
\begin{figure}[h!]
    \centering
    \begin{subfigure}[b]{0.45\textwidth}
        \centering
        \includegraphics[height=0.28\textheight]{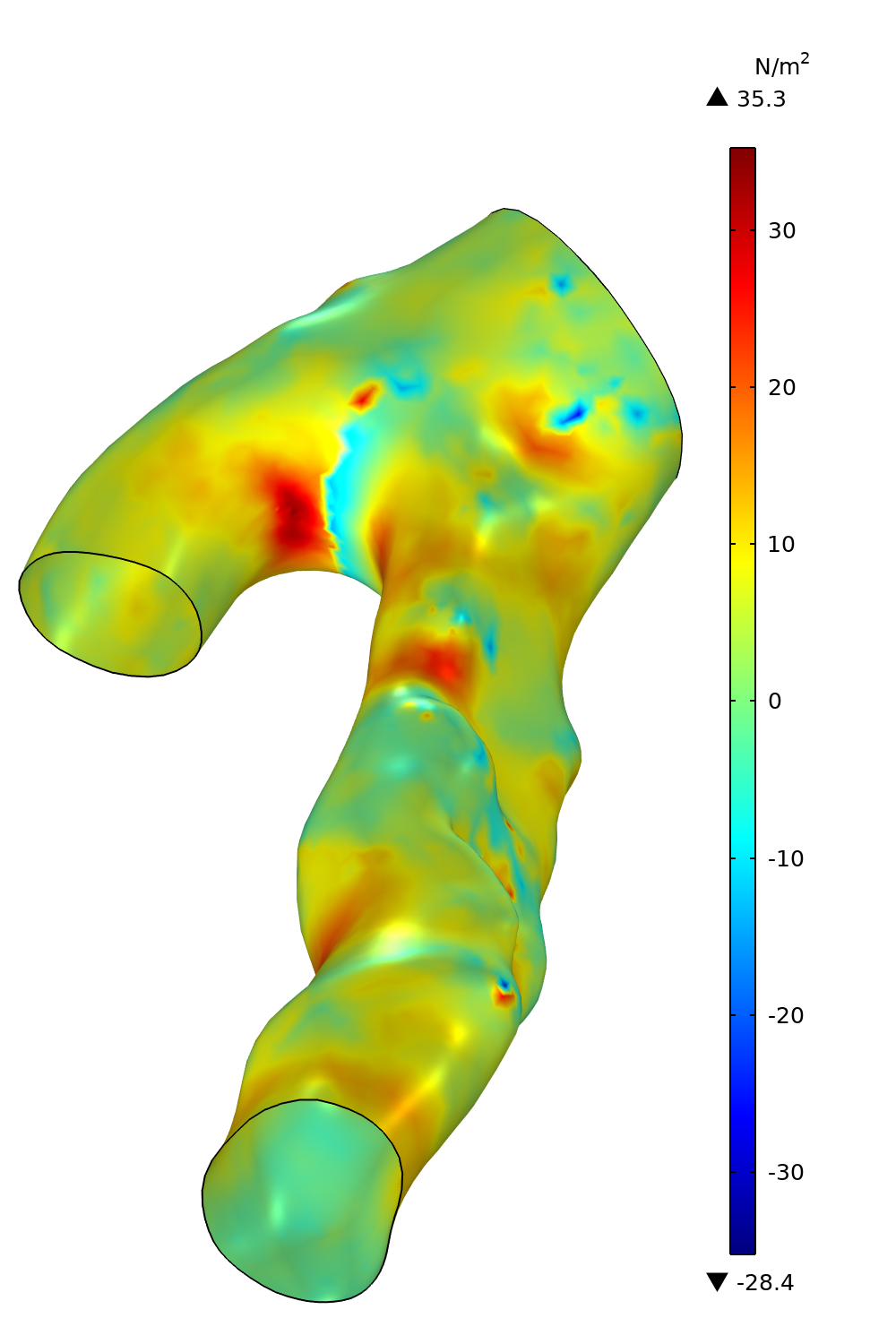}\\ 
        \caption{}
        \label{fig:WSS(a)}
    \end{subfigure}
    \begin{subfigure}[b]{0.45\textwidth}
        \centering
        \includegraphics[height=0.25\textheight]{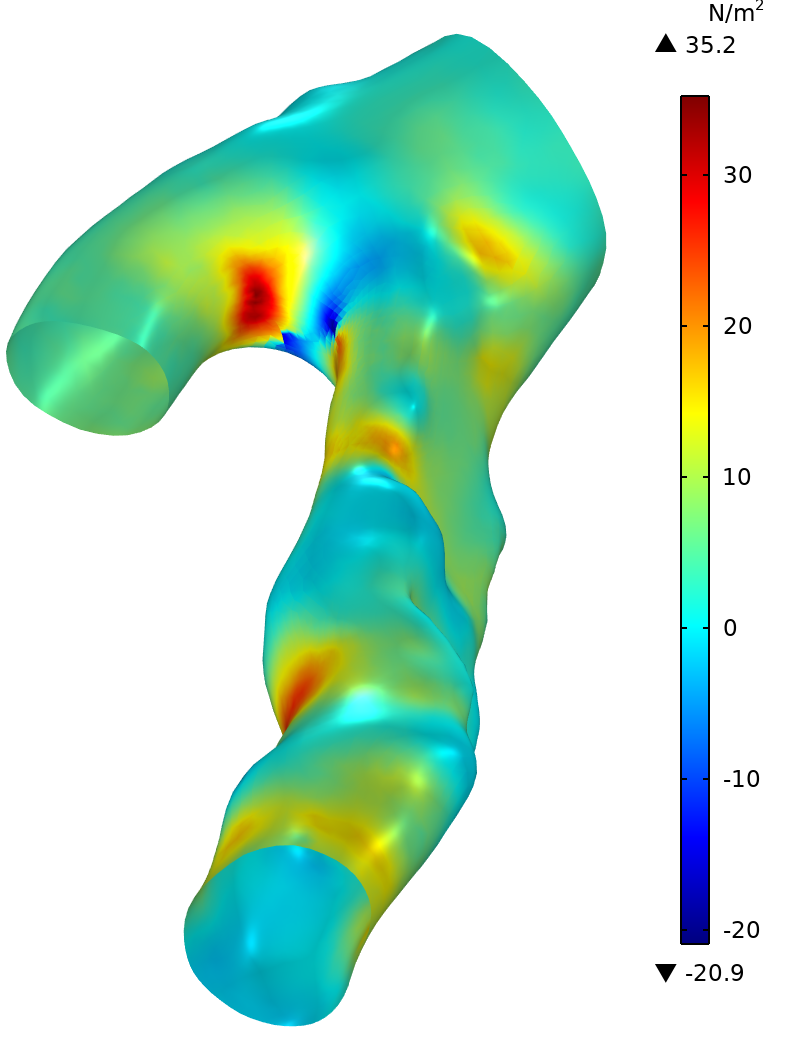}
        \caption{}
        \label{fig:WSS(b)}
    \end{subfigure}
    \begin{subfigure}[b]{0.45\textwidth}
        \centering
        \includegraphics[height=0.25\textheight]{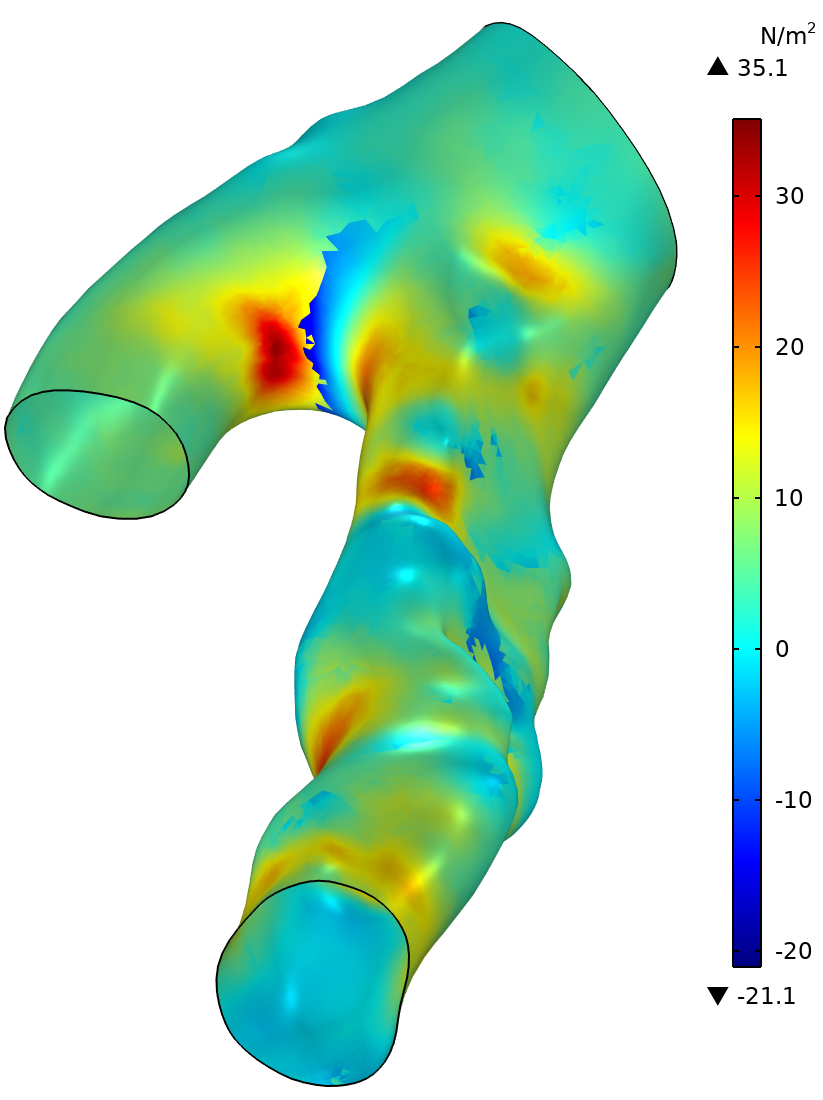}
        \caption{}
        \label{fig:WSS(c)}
    \end{subfigure}
    \begin{subfigure}[b]{0.45\textwidth}
        \centering
        \includegraphics[height=0.25\textheight]{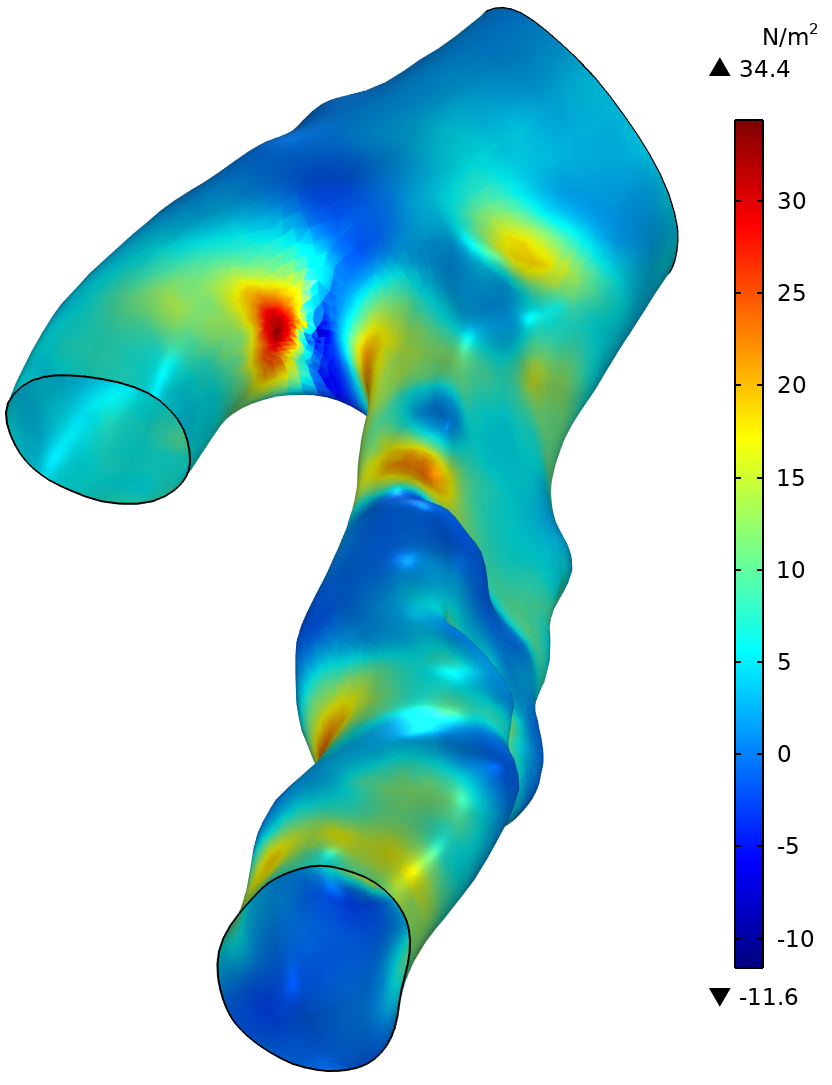}
        \caption{}
        \label{fig:WSS(d)}
    \end{subfigure}
    \caption{Longitudinal WSS distribution  at t=1.1s for all configurations, frontal viewpoint.
    }
    \label{fig:wss_surface}
\end{figure}

\begin{figure}[h!]
    \centering
    \begin{subfigure}{0.45\textwidth}
        \centering
        \includegraphics[height=0.35\textheight]{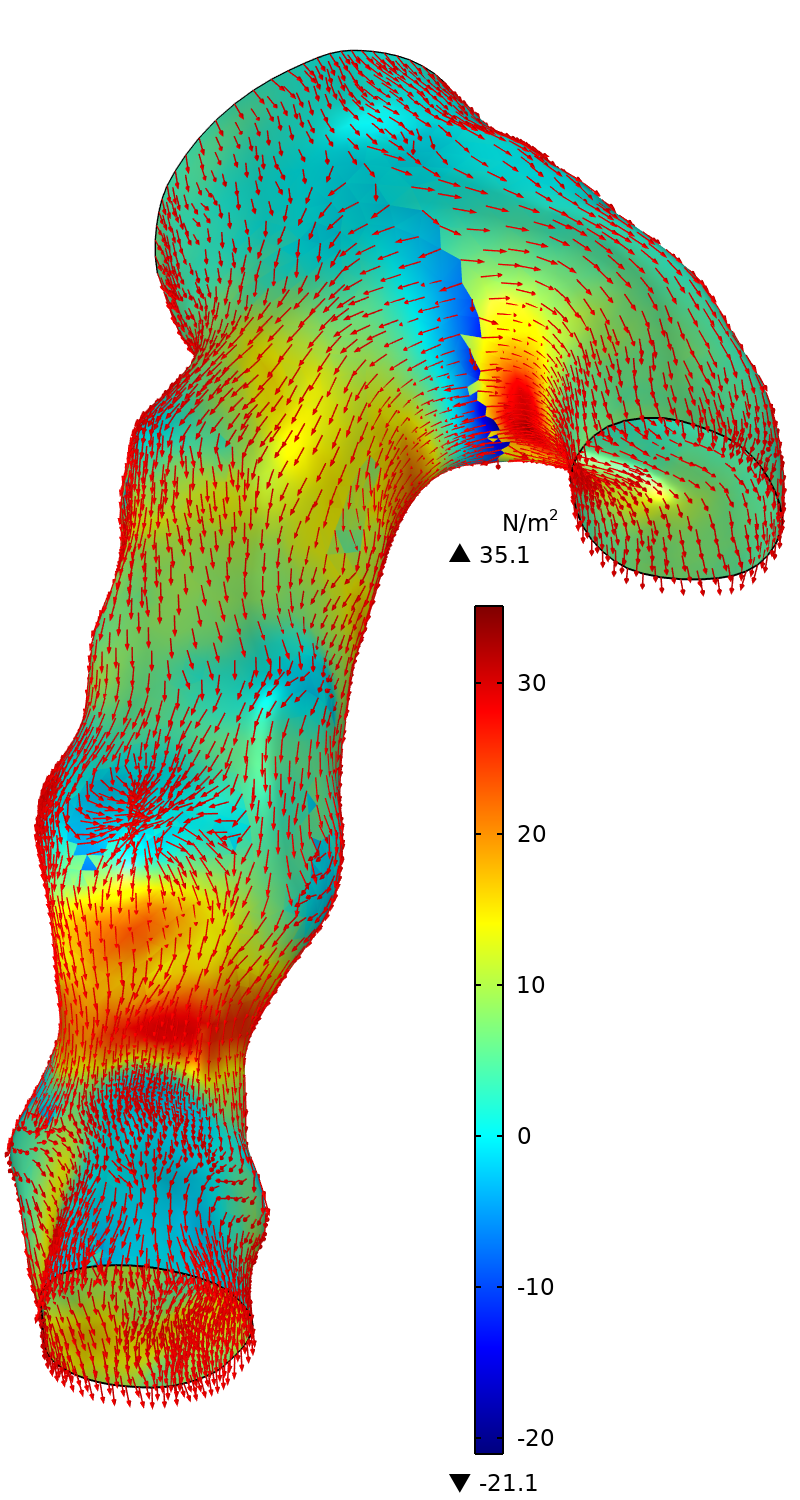}
       \caption*{(c)}
    \end{subfigure}
    \begin{subfigure}{0.45\textwidth}
        \centering
        \includegraphics[height=0.35\textheight]{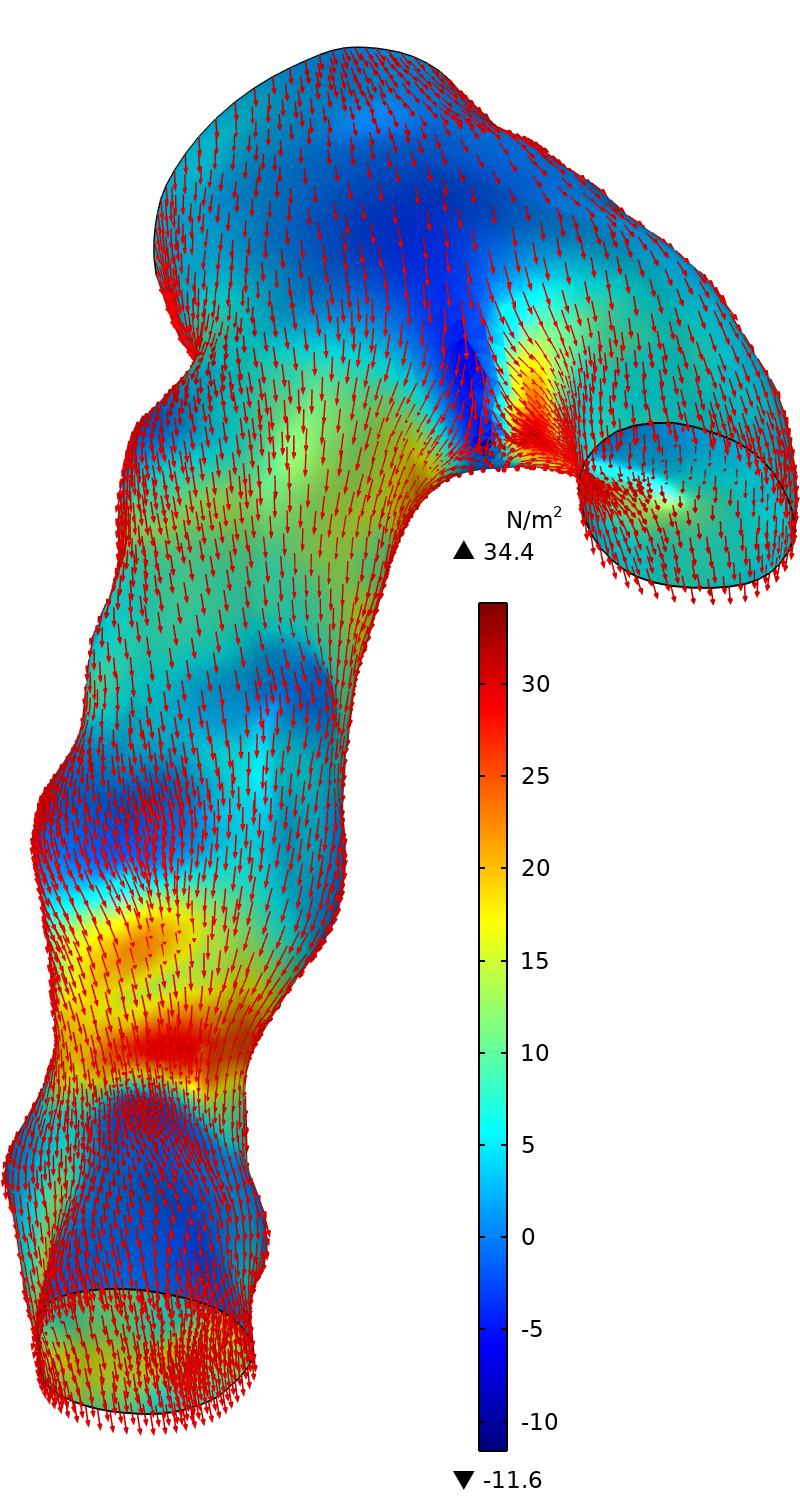}
      \caption*{(d)}
        \end{subfigure}
    \caption{Surface distribution of longitudinal WSS  evaluated using flipped (c) and projected (d) tangent fields (shown as normed arrows),  evaluation with compliant walls at $t=1.1$ s, viewpoint from back (different viewpoint as in Fig. \ref{fig:wss_surface}).}
    \label{fig:flippedtangents_bifurcation}
\end{figure}

\noindent
In Fig. \ref{fig:wss_surface} the surface distributions of $\tau_w^\ell$ for the four configurations at the time of highest flow-rate are presented. 
All plots {(a)}-(d) show high positive WSS up to 35 $\frac{N}{m^2}$ in the bifurcation and in the sinusoidal constrictions of the stenotic bulge of the ICA. 
The center of the bifurcation and the regions around the stenotic bulges of the ICA are regions of low and negative WSS. 
In configuration (a), local point-wise extreme values of $\tau_w^\ell$ appear using the flipped tangent vectors on rigid surface. These are smoothed out in (c) on the deformed surface with elastic walls. Beside these very local phenomena in (a), 
wider areas of negative extreme values occur close to  the separation point of the bifurcation in configurations (b) and (c), compared to (d).
The occurrence of these extreme values is associated with the alignment of chosen tangential fields, which differ for the considered model configurations at the bifurcation point.  
We demonstrate this coherence for automatically rendered flipped (c) and projected tangents (d) on compliant walls in Fig. \ref{fig:flippedtangents_bifurcation}. 
Indeed,  a clear side separation of the  flipped tangent vectors loosing their longitudinal alignment even before the bifurcation can be observed in plot (c),  explaining the discontinuity and  the appearance  of negative WSS values up to -20 $\frac{N}{m^2}$ in this area. 
On the other hand, the longitudinal continuance of  projected tangents  until  the separation point is apparent in plot (d).
Obviously, configurations with flipped tangent vectors lead to a spurious longitudinal WSS evaluation
on the carotid surface close to their bifurcation point, 
whereas the tangent field (d), which is projected from the centerline, seems to appropriately map the main flow and its separation in this problematic area.

\begin{figure}[h!]
    \centering
     \includegraphics[height=0.2\textheight]{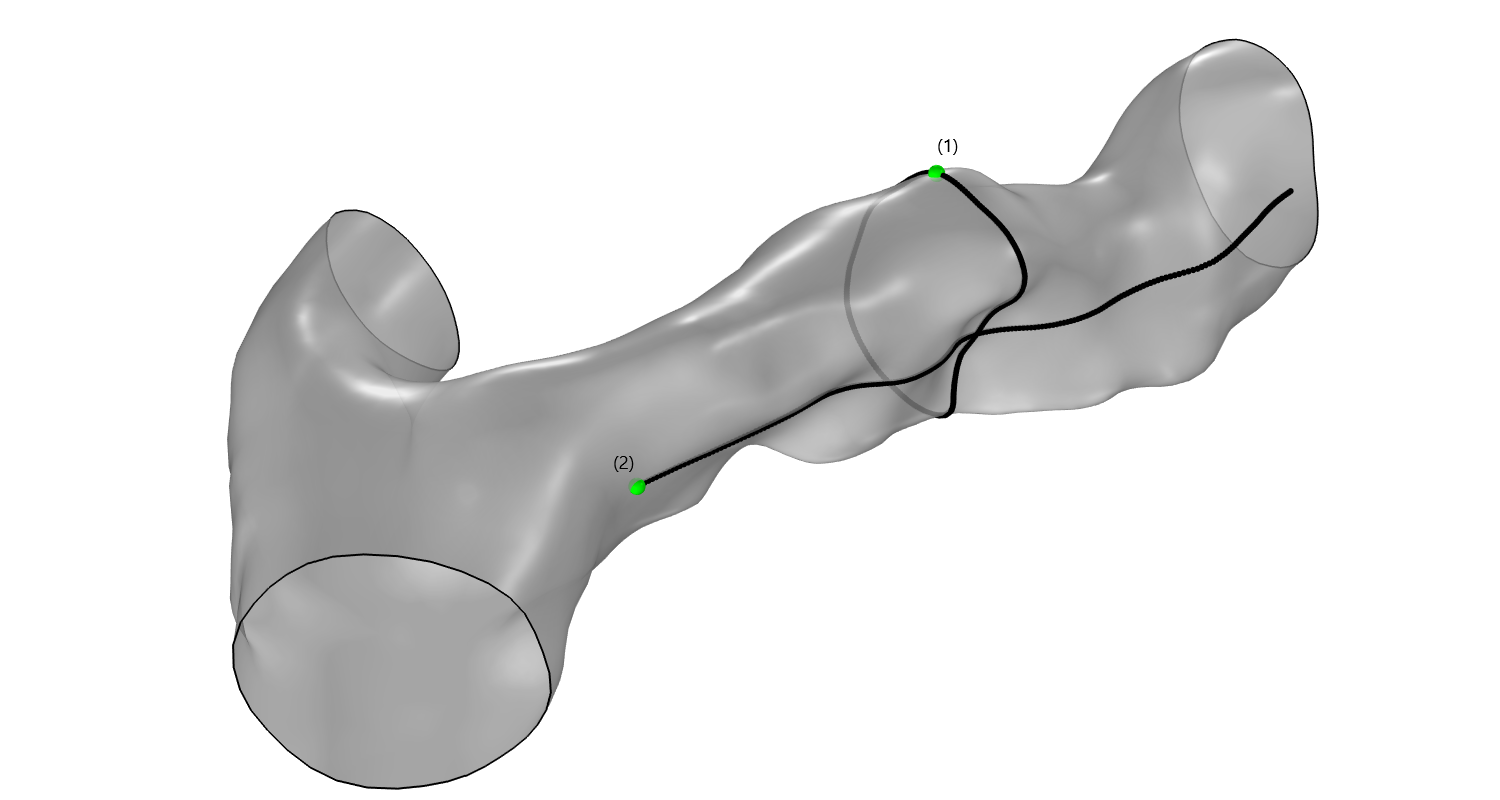}\\
\caption*{Surface intersection curves (1) and (2), frontal viewpoint} \vspace{3mm}
    \begin{subfigure}{0.49\textwidth}
         \includegraphics[height=0.2\textheight,width=\textwidth]{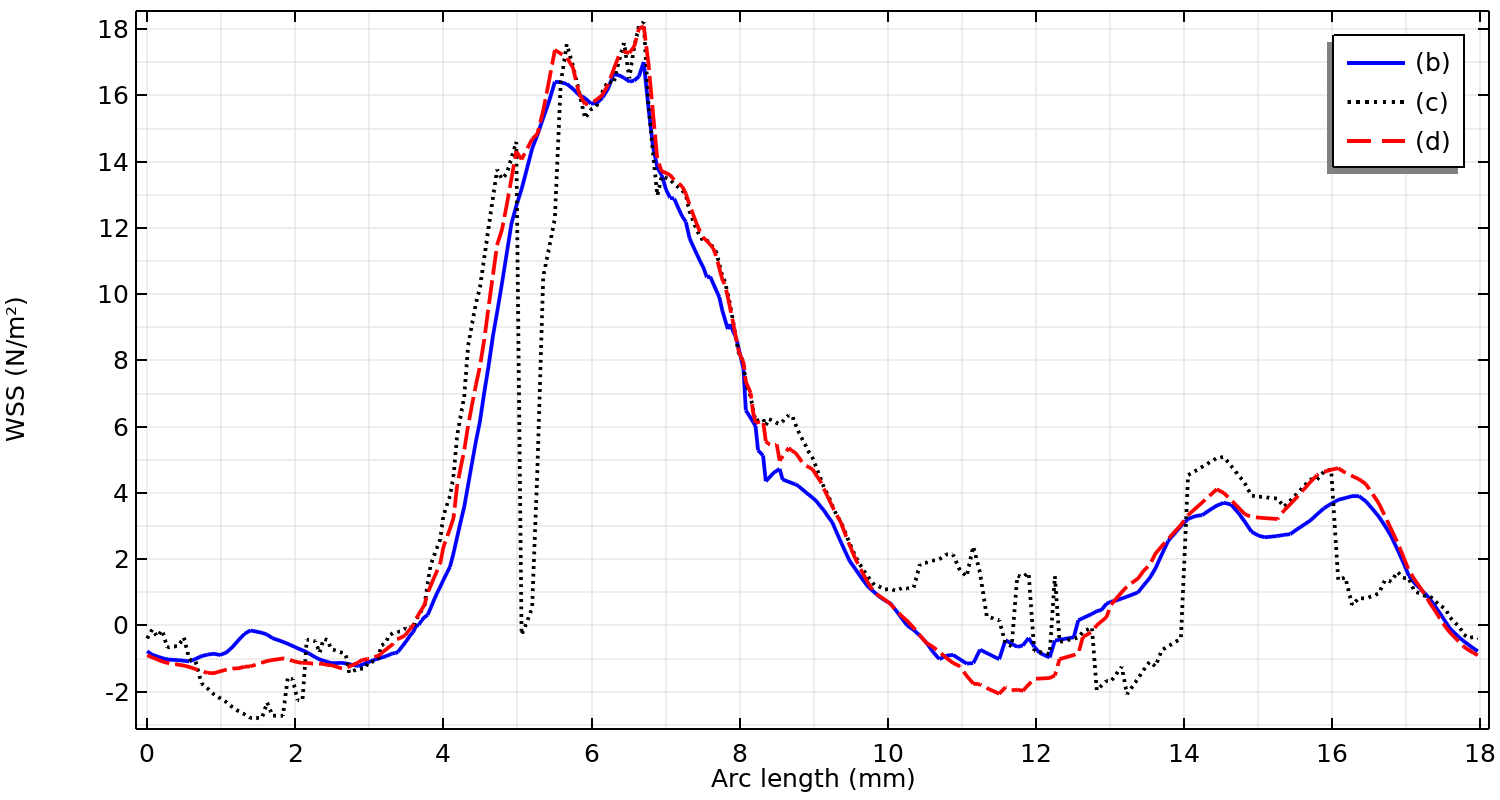}
        \caption*{Intersection  curve (1)}
        \label{fig:WSS_line(1)}
    \end{subfigure}
    \vspace{0.2cm}
    \begin{subfigure}{0.47\textwidth}
        \includegraphics[height=0.2\textheight,width=\textwidth]{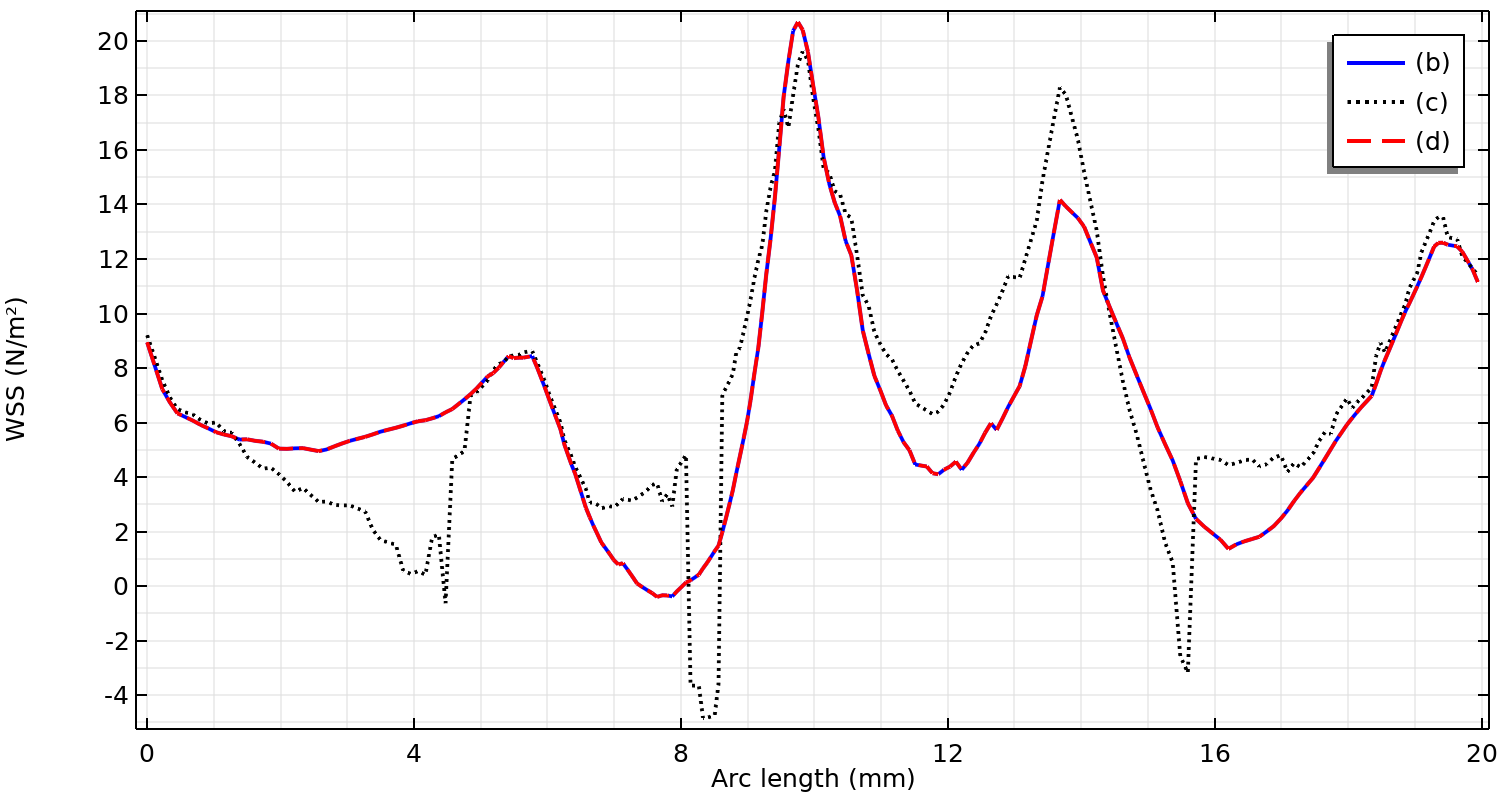}
        \caption*{Longitudinal intersection  curve  (2) }
        \label{fig:WSS_line(2)}
    \end{subfigure}
    \caption{ 
    Upper: circumferential (1) and  longitudinal (2) intersection curves for the evaluation of WSS. The arc length of (1) starts at the respective green point and follow a clockwise direction.
    Below: comparison of  longitudinal WSS for configurations (b)-(d)  along intersection  curves (1),(2). }
    \label{fig:wss_line}
\end{figure}

To get more detailed comparisons, $\tau_w^\ell$ along the chosen surface curves is shown in Fig. \ref{fig:wss_line}  for configurations (b), (c), (d). 
Its values in configurations (b) and (d) are almost identical and differ at most by $1 \frac{N}{m^2}$ along the circumferential curve (1). On the longitudinal line (2) we observe a very good agreement of (b) and (d). The situation is slightly different in case (c) with flipped tangents, 
 where local deflections of WSS values can be observed. These are caused by the previously discussed alignment of differently flipped tangent vectors presented in Fig. \ref{fig:flippedtangents_bifurcation}. The negative deflections of $\tau_w^\ell$ along the longitudinal curve (2) are not visible in Fig. \ref{fig:flippedtangents_bifurcation} due to a different viewpoint, nevertheless they can be identified with  local blue areas on the outer wall of ICA in plot (c) of Fig. \ref{fig:wss_surface}.

Generally, we can conclude that the choice of the tangential field has a considerable effect on the evaluation of longitudinal WSS  in the bifurcation region of the carotid artery surface, with 
the projected tangent vectors being the most suitable for this analysis. 
After passing the problematic area of the bifurcation the WSS results for projected and flipped tangent fields are comparable at maximal systolic flow, whereas certain inaccuracies occur using automatic flipped tangent vectors. This is because of the misalignment in longitudinal direction on complex and uneven surfaces.
In addition, effects of the wall movement in FSI models, are present around the bifurcation point, compare e.g., plot (b) for rigid and (d) for compliant walls in Fig. \ref{fig:wss_surface}.
 
\subsubsection{Vector-valued WSS}

\begin{figure}[h!]
    \centering
    \begin{subfigure}{0.45\textwidth}
    \begin{overpic}[scale=.53,,tics=10]
    {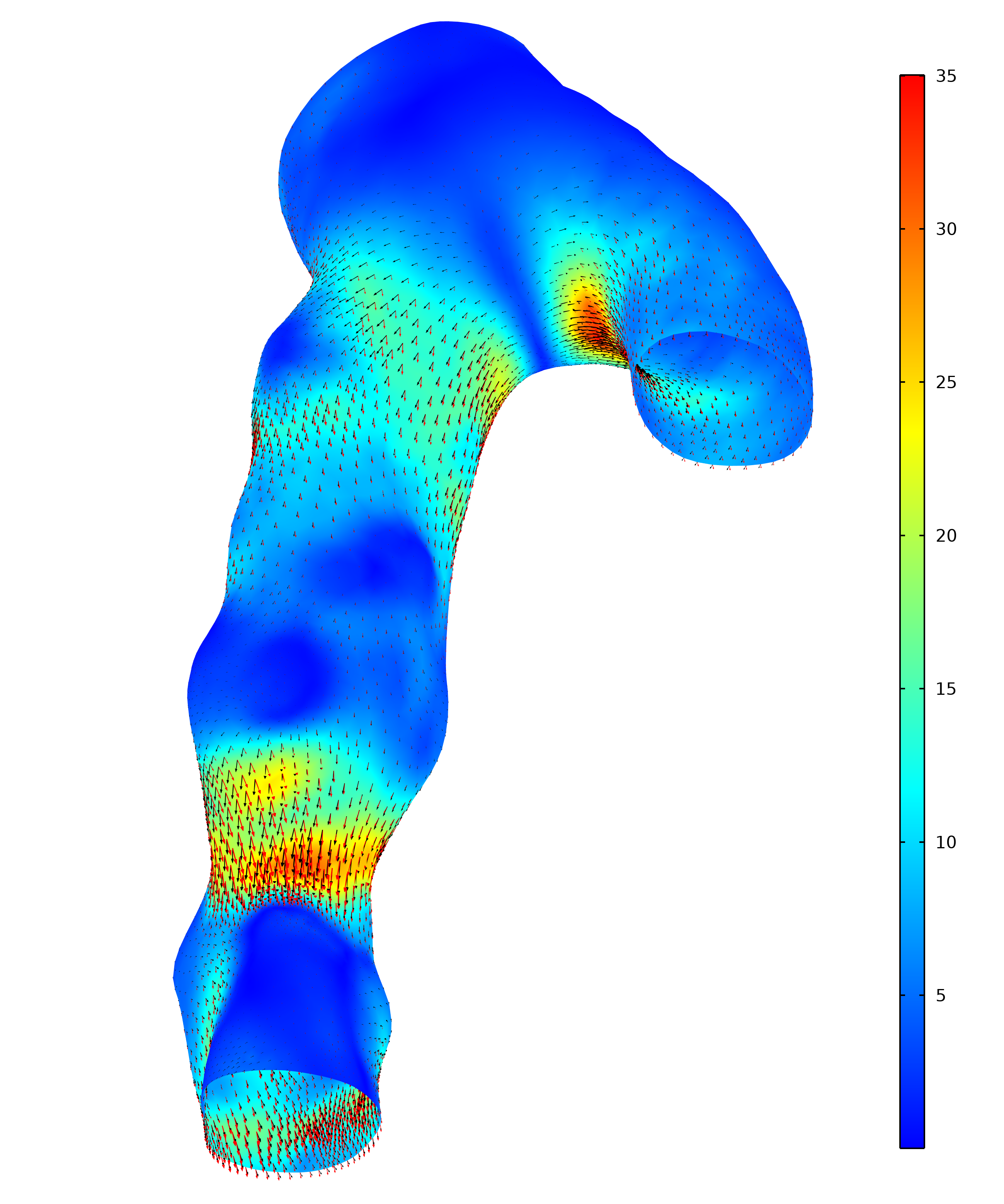}
        \put(28,56){\linethickness{0.25mm}\color{black}\polygon(0,0)(25,0)(25,27)(0,27)}
    \end{overpic}
    \end{subfigure} \qquad  \qquad 
    \begin{subfigure}{0.35\textwidth}
         \includegraphics[height=0.25\textheight]{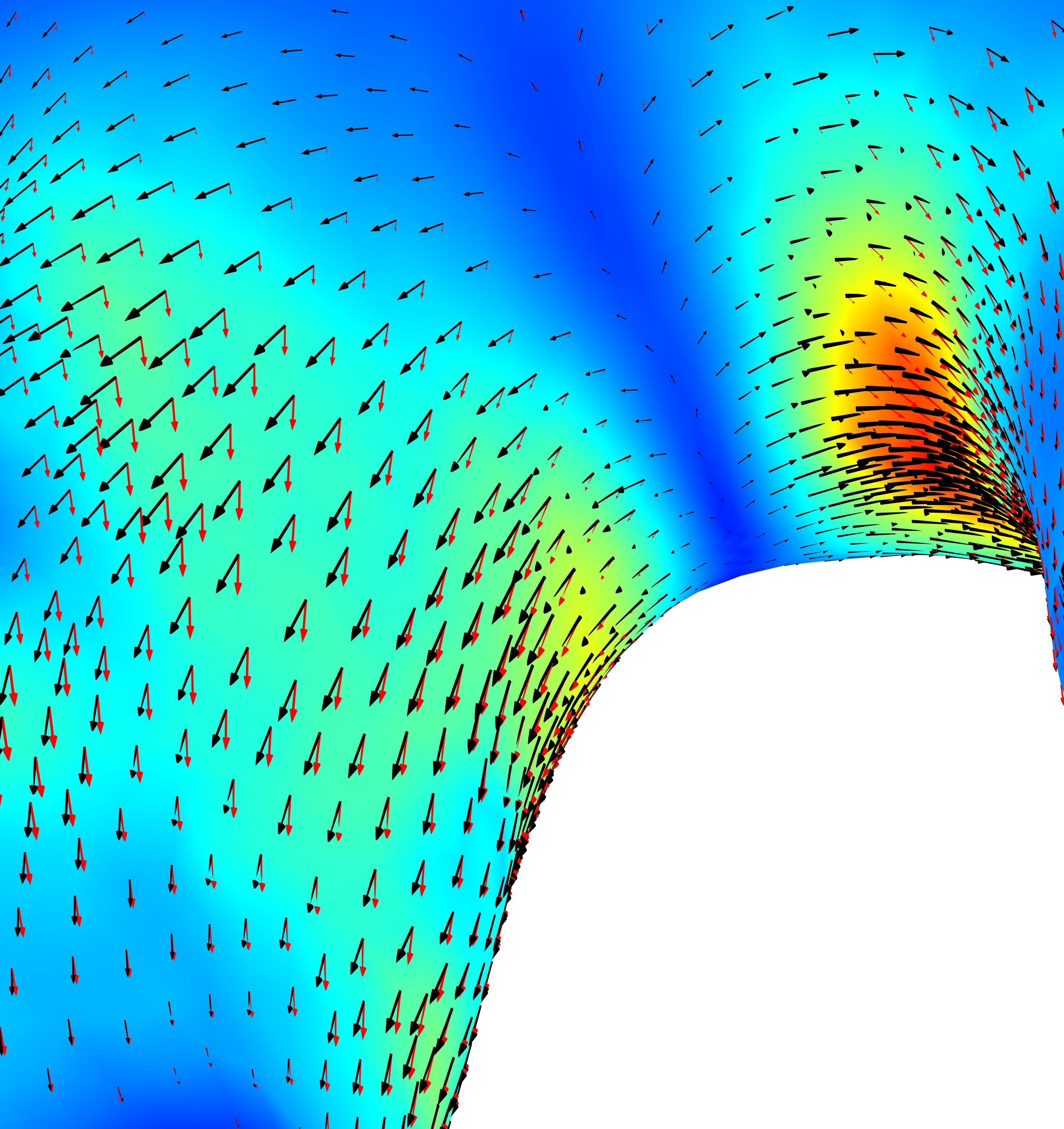}
        \caption*{Detail of the bifurcation point: WSS vector $\vec{\tau}_w$  (black arrows) and its longitudinal component $\tau_w^\ell$ (red arrows).}
    \end{subfigure}
    \caption{Surface: $\tau_w^a$ - amplitude of the WSS vector on compliant walls at $t=1.1$ s. Arrows: WSS vector   $\vec{\tau}_w$ (black) and its longitudinal component  $\tau_w^\ell$ (red), both proportional to its size.}
    \label{fig:WSS_vector}
\end{figure}

The WSS on the carotid inner surface measured by the amplitude of the WSS vector $\tau_w^a$   (\ref{wss1}) is presented in Fig. \ref{fig:WSS_vector} with  corresponding arrows of  vector $\vec{\tau}_w$  as well as its longitudinal components ${\tau}_w^\ell$ at the time point of maximal flow in the reference geometry frame. One can observe similar surface distribution of the amplitude of WSS  and of the size of its longitudinal component comparing the Figures \ref{fig:WSS_vector} and \ref{fig:flippedtangents_bifurcation} - (d). Let us note, that the opposite direction of $\vec{\tau}_w$ with respect to the main flow, which occurs sometimes in blue areas with low WSS amplitude, cannot be recognised by $\tau_w^a$. Thus the backward flow cannot be detected in the surface distribution on Fig.  \ref{fig:WSS_vector}. 
The negative sign and the sign change of WSS resulting from opposite alignment with respect to the main flow plays a role in the evaluation of the oscillatory flow behavior, and it is discussed in what follows.

\subsection{Oscillatory shear index}

In this section we present the temporal change of both, the longitudinal WSS and the vector-valued WSS during the whole cardiac cycle in means of the oscillatory shear index and compare them for compliant as well as rigid wall models. The results are presented on the inner arterial surface in the reference, i.e., the initial geometry frame.

\begin{figure}[h!]
    \centering
    \begin{subfigure}{0.2\textwidth}
        \includegraphics[height=0.35\textheight]{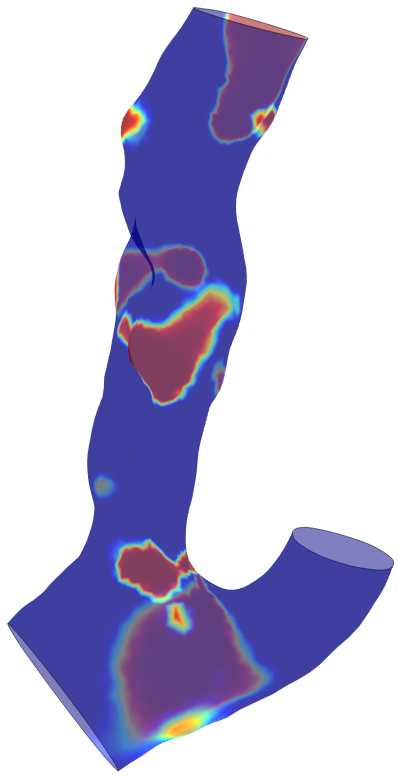}
        \caption*{(b)}
        \label{fig:osi(b)}
    \end{subfigure}
    \hspace{0.06\textwidth}
    \begin{subfigure}{0.22\textwidth}
        \includegraphics[trim=0 0 220 0, clip,height=0.35\textheight]{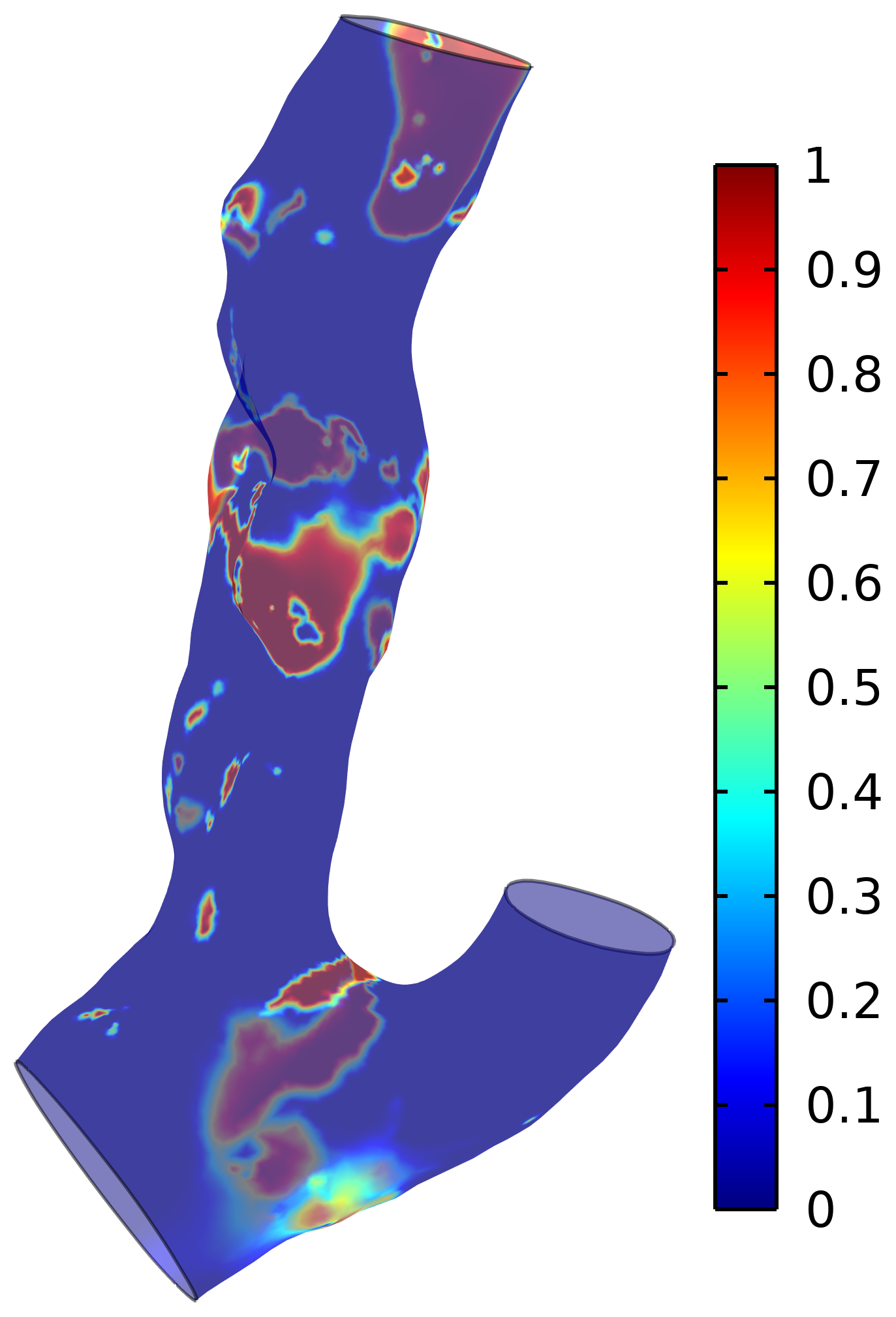}
        \caption*{(c) }
        \label{fig:osi(c)}
    \end{subfigure}
    \hspace{0.06\textwidth}
    \begin{subfigure}{0.22\textwidth}
        \centering
        \includegraphics[height=0.35\textheight]{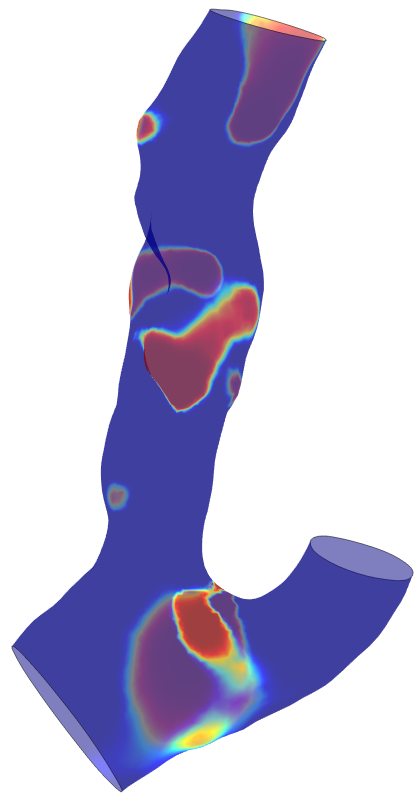}
        \caption*{(d)}
        \label{fig:osi(d)}
    \end{subfigure}
    \hspace{0.06\textwidth}
   \begin{subfigure}{0.1\textwidth}
        \includegraphics[trim=800 0 0 0, clip,height=0.3\textheight]{Bilder/OSI__c_.png}
        \caption*{ }
    \end{subfigure}
    \caption{Oscillatory shear index  $OSI^\ell$ (\ref{osi1}) for longitudinal WSS  ($\tau_w^\ell$)  evaluated for the second cardiac cycle, $t\in [0.9s,1.8s]$ in model configurations (b),(c),(d) see Table \ref{table:configurations}.   Red areas show the occurrence of long lasting negative $\tau_w^\ell$.}
    \label{fig:OSI}
\end{figure}

At first, the oscillatory 
index (\ref{osi1}) 
for the longitudinal WSS (\ref{wss2})  is presented  for the fixed wall model using projected tangents (b) and for the FSI model using flipped (c) and projected tangents (d) in Fig. \ref{fig:OSI}.
Comparing the results, 
almost no difference in the OSI evaluation for configurations (b) and (d) with projected tangents can be observed.  The only region of difference worth mentioning spreads out in the bifurcation area, where  the centerline tangents have been projected on different surfaces, obviously  due to the wall deformation in case (d).
In contrast to (b) and (d), configuration (c) shows many  small-scale and point-like 
 extreme values, which is 
 a consequence of the discontinuity  and deflections  of WSS  arising from the erratic alignment of the flipped tangents in some surface regions discussed above.

 Concerning the hemodynamical interpretation of the results presented in Fig. \ref{fig:OSI},  conspicuous red regions of maxima indicating long-lasting negative WSS inside and oscillating WSS at their edges (green transition zones) can be observed in the bifurcation zone as well as prior and posterior to the sinusoidal stenotic occlusion of the ICA in all model configurations. Further punctual abnormalities are present, e.g., after the stenotic occlusion on the left side of the ICA wall. 
 In consistency with the streamlines presented in Fig. \ref{fig:streamlines}, red $OSI^\ell$ regions are related to vortices adjacent to the stenosis bulges of the ICA or prior to the bifurcation.  
 Those maximum regions are an indicator of static reversal flow (vortices) and the green transition zones of high longitudinal WSS oscillations  indicate  the pathological progression of mechanical damage of the  artery wall, according to the hemodynamic hypothesis.
 
\begin{figure}[h!]
    \centering
     \includegraphics[height=0.45\textheight]{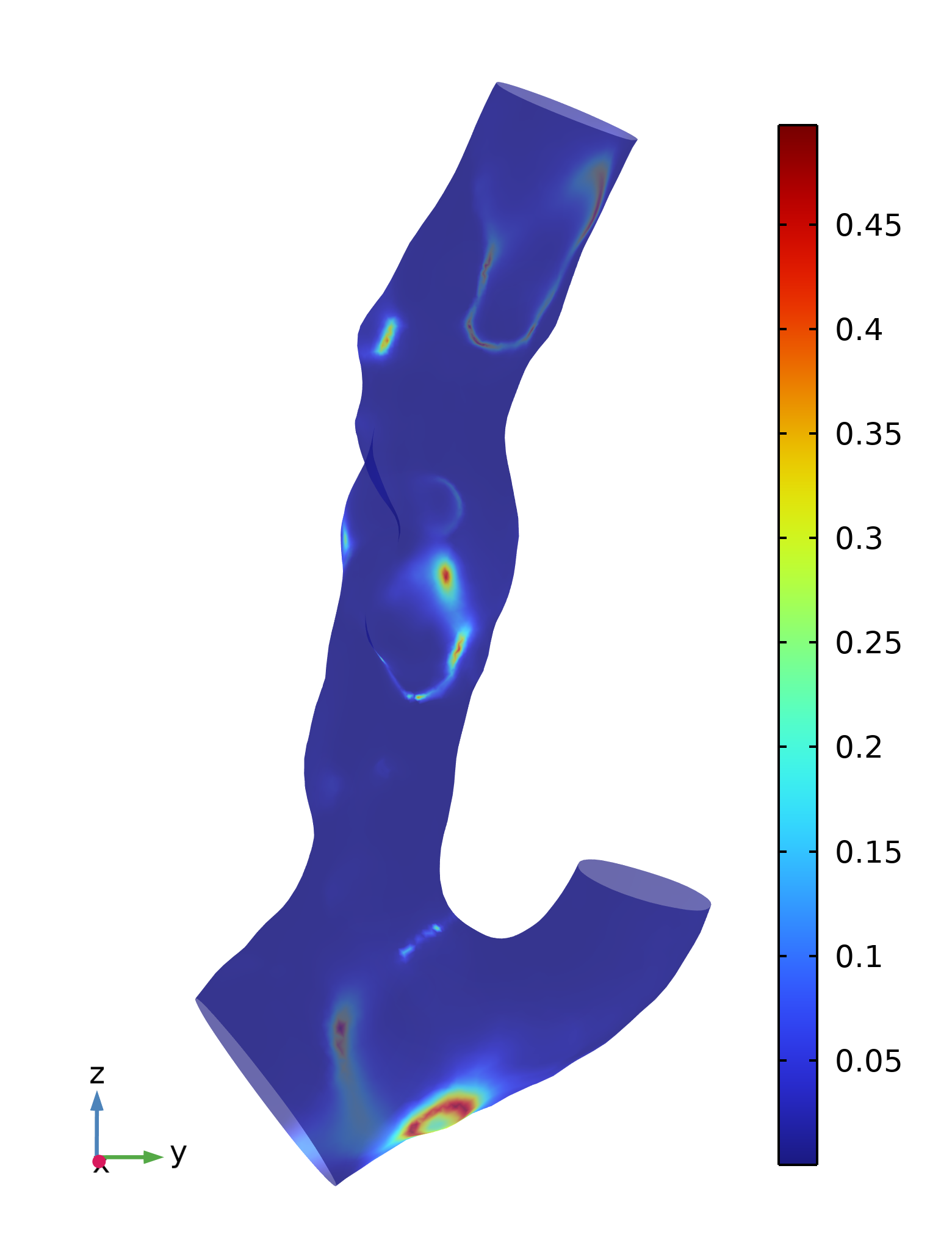}
        \caption{$OSI$ index (\ref{osi2})  for vector-valued WSS ($\vec{\tau_w}$)  evaluated for the second cardiac cycle ($t\in [0.9s,1.8s]$), compliant wall model.}
        \label{fig:osi2}
    \end{figure}

 Finally, a second oscillatory shear index (\ref{osi2}) for the vector-valued  WSS computed using  the compliant wall  model is presented in Fig. \ref{fig:osi2}. Here, the red maximum values correspond to WSS regions with zero temporal mean and represent zones 
 of complete sign balance of  $\vec{\tau}_w$, indicating temporal oscillations. As mentioned in Section \ref{subsec:osi}, static vortices and continuous recirculations, corresponding to a permanent negative sign of the WSS, are zero valued in this OSI definition and cannot be tracked here. Nevertheless, the high OSI 
 edge-like regions in Fig. \ref{fig:osi2} are in good consistency with the green transition zones observed in Fig. \ref{fig:OSI} - (b),(d), both indicating  high wall shear stress oscillations  and the pathological progression e.g. on the right and back side of the common carotid artery prior to the bifurcation point, or on the upper back side of the considered ICA region.
 Small discrepancies between the maxima edges of $OSI$ and the green transition zones of $OSI^\ell$ in Figs. \ref{fig:OSI} and  \ref{fig:osi2} are obviously caused by differences in the considered wall shear stress indicators $\vec{\tau}_w$ and ${\tau}_w^\ell$.
 
 \section{Conclusion}

In this contribution a computational study of fluid dynamic and hemodynamic risk parameters has been performed for a carotid artery. The patient-based lumen and its surrounding walls have been imported into the numerical software {\sl Comsol Multiphysics} and were supplemented for one with fluid dynamics confined by rigid walls as well as a FSI model for compliant artery using the
linear elasticity deformation model. 
The  effects of fluid stresses on the arterial wall have been quantified using established hemodynamic risk factors: wall shear stress and oscillatory shear index. 
Following the low shear theory, we focused on exploration of low and negative WSS and evaluated the longitudinal component  of WSS vector (longitudinal WSS) allowing to track the reverse flow in patient-based morphology.
The presented results demonstrate the strong dependency of the orientation-based longitudinal WSS and its OSI-index on the proper construction 
of tangential vectors and address this problem on topologically complex surfaces obtained from patient CTA scans. 
For the studied carotid artery tree, we applied the projection of centerline tangents to the inner arterial surface in order to obtain  properly aligned and  smooth tangential field, and compared the longitudinal WSS computed for projected as well as generic mesh-based tangent vectors.
Since the projected tangential field retains the longitudinal alignment on the craggy surface and maps the flow separation in the bifurcation area much better than the automatically generated tangent vectors, reliable numerical results for longitudinal WSS, allowing hemodynamic predictions of reverse flow have been obtained by applying the projection method.  

For comparisons, the commonly used vector-valued WSS and its amplitude have been evaluated as well and its oscillatory behavior  has been quantified by the corresponding OSI index.
Since the WSS-vector amplitude and the amplitude of its temporal mean does not track its opposing orientation,  reversal flow or vortices cannot be detected, which is confirmed by our numerical results. Even though the OSI index, based on the vector-valued WSS instead of its longitudinal component, indicates flow oscillatory regions well, it is not able to detect long lasting reversal flow regions with low and negative WSS.
Hence, the major benefit of our modeling approach using the direction-based longitudinal WSS and its oscillatory index lies the investigation of sites of low shear stresses and persistent reversal flow.

Obtained numerical data and derived risk indicators can be complemented with measured medical patient data and they will be explored using machine learning algorithm within the joint research project MLgSA, see Acknowledgement, to assess the potential stroke risk for a patient.

For further work, the range of risk prediction tools can be extended by considering further multi-directional WSS parameters, where the choice of the tangential field is of particular importance. According to
\cite{hoogendoorn},  time-averaged WSS (TAWSS; \cite{arzani, john}) and  relative residence time (RRT; \cite{gorring, hashemi, soulis}) are strong predictive clinical markers for disease development, even in early stages of atherosclerosis. 
Moreover, transversal WSS (\cite{mohamied16, pedrigi, pfeiffer15}) and cross flow index (CFI, \cite{mohamied16}) have shown a predictive value when complex fluid flow appears in later phases of atherosclerosis. Taking those parameters into account, a more differentiated prediction of atherosclerotic development can be achieved.

\section*{Acknowledgements}
This work was funded by the  German Federal Ministry of Education and Research (BMBF) joint Project \href{https://math4innovation.de/index.php?id=70}{05M20UNA-MLgSA}.


\subsection*{Literature}

\begin{thebibliography}{20}
 
 \bibitem{antiga}
		Antiga, L.:Patient-specific Modeling of Geometry and Blood Flow in Large Arteries. {\it PhD thesis}, Politecnico di Milano, 2002
		
		
		\bibitem{arzani} Arzani, A., Shadden S.C.:
		Characterizations and correlations of wall shear stress in aneurysmal flow,
		{\it Journal of Biomechanical Engineering}, 138, 0145031-01450310, 2016,
		DOI: \href{https://doi.org/10.1115/1.4032056}{10.1115/1.4032056}
		
		\bibitem{blagojevic}
		Blagojević, M., Nicolić, A., Zivković, M., Zivković, M., Stanković, G.,Pavlović, A.:
		Role of oscillatory shear index in predicting the occurrence and development of plaque,
		{\it Journal of the Serbian Society of Computational Mechanics}, 7(2), 29-37, 2013
		
		\bibitem{bukacsplitting}
		 Buka\v c, M.,  \v Cani\'c, S.,  Glowinski, R.,  Muha B.,  Quaini, A.:: A modular, operator-splitting scheme for fluid–structure
		interaction problems with thick structures, {\it International journal for
			numerical methods in fluids} 74.8, 577–604, 2014, DOI: \href{https://doi.org/10.1002/fld.3863}{10.1002/fld.3863}
		
		
		\bibitem{canicshell}
		\v Cani\'c, S., Tamba\v ca, J., Guidoboni, G., Mikeli\'c, A., Hartley, C.J., Rosenstrauch, D.:
		Modeling viscoelastic of arterial walls and their interaction with pulsatile blood flow,
		{\it SIAM J. Appl. Math.} 67(1), 164–193, 2006
		
		\bibitem{ciarletshell}
		Ciarlet, P.G.: {\it Mathematical Elasticity, Volume III: Theory of Shells}, North-Holland,
		Amsterdam, 2000
		
		
		\bibitem{comsol} Comsol Multiphysics. Reference Manual,  
		\\ \href{https://doc.comsol.com/5.6/doc/com.comsol.help.comsol/COMSOL_ReferenceManual.pdf}{https://doc.comsol.com/5.6/doc/com.comsol.help.comsol/COMSOL-ReferenceManual.pdf},  2022 
		
		\bibitem{debus} Debus, E.S., Torsello, G., Schmitz-Rixen, G., Flessenkämper, I., Storck, M., Wehk. H., Grundmann, R.T.:
		Ursachen und Risikofaktoren der Arteriosklerose,
		{\it Gefässchirurgie}, 18, 2413-2419, 2013,
		DOI: \href{https://doi.org/10.1007/s00772-013-1233-6}{10.1007/s00772-013-1233-6}
		
		
		\bibitem{Eulzer1}
		Eulzer, P., Meuschke, M., Klinger,
		C. M., Lawonn, K.: Visualizing carotid blood flow
		simulations for stroke prevention, {\it Computer Graphics Forum}, 40
		435–446, 2021, DOI:
		\href{https://doi.org/10.1111/cgf.14319}{10.1111/cgf.14319}
		
		\bibitem{Eulzer2} Eulzer, P., Richter, K.,  Meuschke M.,  Hundertmark A.,  Lawonn K.: Automatic Cutting and Flattening of Carotid Artery Geometries, {\it in  Proceedings of Eurographics Workshop on Visual Computing for Biology and Medicine 2021}, (Editors: S. Oeltze-Jafra, N. N. Smit, and B. Sommer)
		
		\bibitem{gallo}Gallo, D., Steinman, D.A.,  Morbiducci, U.:
		Insights into the co-localization of magnitude-based versus direction-based indicators of disturbed shear at the carotid bifurcation,
		{\it Journal of Biomechanics}, 49(2), 2413-2419, 2016,
		DOI: \href{https://doi.org/10.1016/j.jbiomech.2016.02.010}{10.1016/j.jbiomech.2016.02.010}
		
		\bibitem{gorring}
		Gorring, N., Kark, L., Simmons, A., Barber, T.:
		Determining possible thrombus sites in an extracorporeal device, using computational fluid dynamics-derived relative residence time.
		{\it  Computer Methods in Biomechanics and Biomedical Engineering}, 18(6), 628-634, 2014,
		DOI: \href{https://doi.org/10.1080/10255842.2013.826655}{10.1080/10255842.2013.826655}
		
		
		\bibitem{hashemi}
		Hashemi, J., Patel B., Chatzizisisi, Y.S., Kassab, G.S.:
		Study of Coronary Atherosclerosis Using Blood Residence Time
		{\it Front. Physiol.}, 12, 625420, 2021, 
		DOI: \href{https://doi.org/10.3389/fphys.2021.625420}{ 10.3389/fphys.2021.625420}
		
		\bibitem{hoogendoorn} Hoogendoorn, A., Hartman, E. MJ., Kok, A., De Nisco, G.:
		Multidirectional wall shear stress promotes advanced coronary plaque development - comparing five shear stress metrics,
		{\it Cardiovascular Research}, 116(6), 1136-1146, 2019,
		DOI: \href{https://doi.org/10.1093/cvr/cvz212}{10.1093/cvr/cvz212}
		
		
		
		\bibitem{Hundertmark2010}  Hundertmark, A.,  Luk\'a\v cov\'a, M.: Numerical study of shear-dependent non-Newtonian fluids in compliant vessels, {\it Computers and Mathematics with Applications}, 60,   572-59, 2010.
		
		\bibitem{hundertmarkshear}
		Hundertmark, A., Luk\'a\v cov\'a, M., Rusn\'akov\'a, G.: Fluid-structure
		interaction for shear-dependent non-Newtonian fluids, {\it In: Kaplick\'y P. (Ed.) Topics
			in Mathematical Modeling and Analysis} Vol.7, 109–158, 2012
		
		\bibitem{iannuzzi}Iannuzzi, A., Rubba, P., Gentile, M., Mallardo, V., Calcaterra, I., Bresciani, A., Covetti, G., Cuomo, G., Merone, P., Di Lorenzo, A., Alfieri, R., Aliberti, E., Giallauria, F., Di Minno, M., Iannuzzo, G.:
		Carotid Atherosclerosis, Ultrasound and Lipoproteins.
		{\it Biomedicines}, 9(5), 521, 2021,
		DOI: \href{https://doi.org/10.3390/biomedicines9050521}{10.3390/biomedicines9050521}
		
		\bibitem{izzo2018}
		Izzo, R., Steinman, D., Manini, S., Antiga L.:
		The Vascular Modeling Toolkit: A Python Library for the Analysis of Tubular Structures in Medical Images
		{\it The Open Journal} 25(3), 2018
		
		\bibitem{janelashear}
		Janela, J., Moura, A., Sequeira, A.: A 3D non-Newtonian fluid-structure interaction
		model for blood flow in arteries, {\it J. Comput. Appl. Math. 234}, 2783–2791, 2010
		
		
		\bibitem{john}
		John, L., Pustejovská, P., Steinbach, O.:
		On the influence of the wall shear stress vector form on hemodynamic indicators,
		{\it Computing and Visualization in Science}, 18, 113-122, 2017
		DOI: \href{https://doi.org/10.1007/s00791-017-0277-7}{10.1007/s00791-017-0277-7}
		
		\bibitem{ku} Ku, D. N., Giddens, D. P.,  Zarins,  C. K., . Glagov,  S.:
		Pulsatile flow and atherosclerosis in the human carotid bifurcation, positive correlation between plaque location and low oscillating shear stress. {\it Arteriosclerosis}, 5, 293–302, 1985
		
		\bibitem{Giorgio_low_WSS} 
		Lawton, M., Higashida, R.,  Smith, W. S.,  Young, W. L.,  Saloner, D., Boussel, L., Rayz, V.,  McCulloch, C.,  Martin, A., Acevedo-Bolton, G.,: 
		Aneurysm growth occurs at region of low wall shear stress: Patient-specific correlation of hemodynamics and growth in a longitudinal study, 
		{\it Stroke},  39,  2997-3002,  2008,
		DOI: \href{https://www.ahajournals.org/doi/full/10.1161/STROKEAHA.108.521617}{10.1161/STROKEAHA.108.521617}
		
		\bibitem{mohamied16} Mohamied, Y., Sherwin, S.J., Weinberg, P.D.:
		Understanding the fluid mechanics behind transverse wall shear stress,
		{\it Journal of Biomechanics}, 50(4), 102-109, 2016,
		DOI: \href{https://doi.org/10.1016/j.jbiomech.2016.11.035}{10.1016/j.jbiomech.2016.11.035}
		
		\bibitem{mohamied14} Mohamied, Y., Rowland, E.M., Bailey, E.L. Sherwin, S.J., Schwartz, M. A., Weinberg, P.D.:
		Change of Direction in the Biomechanics of Atherosclerosis,
		{\it Annals of Biomedical Engineering}, 43, 16-25, 2014.
		DOI: \href{https://doi.org/10.1007/s10439-014-1095-4}{10.1007/s10439-014-1095-4}
		
		\bibitem{morbiducci15} 
		Morbiducci, U., Gallo, D., Cristofanelli, S., Ponzini, R., Deriu, M. A., Rizza, G., Steinmann, D. A.: 
		A rational approach to defining principal axes of multidirectional wall shear stress in realistic vascular geometries, with application to the study of the
		influence of helical flow on wall shear stress directionality in aorta, 
		{\it Journal of Biomechanics}, 48(6), 899-906, 2015, 
		DOI: \href{https://doi.org/10.1016/j.jbiomech.2015.02.027}{10.1016/j.jbiomech.2015.02.027}
		
		\bibitem{perktold95}
		Perktold, K., Rappitsch, G.: Computed simulation of local blood flow and vessel mechanics in a compliant carotid artery bifurcation model, {\it Journal of Biomechanics},  28(7), 845–856, 1995
		
		
		\bibitem{Qua_OSI} Quarteroni, A.,  Gianluigi, R.:
		Optimal control and shape optimization of aorto-coronaric bypass anastomoses,
		{\it Math. Models Methods Appl. Sci.}, 13 (12), 1801-1823, 2003
		
		\bibitem{quarteroni04} Quarteroni, A.,  Formaggia, A: Mathematical modelling and numerical 
		simulation of the cardiovascular system.
		{\it In Ciarlet, P.G. $\&$ Lions, J.L. (Hrsg.), Handbook of Numerical Analysis, }12 (3-127), 2004, Amsterdam: Elsevier. 
		
		\bibitem{quarteroni2019}
		Quarteroni, A., Dede', L., Manzoni, A., Vergara, C.:
		Mathematical Modelling of the Human Cardiovascular System: Data, Numerical Approximation, Clinical Applications, {\it Cambridge Monographs on Applied and Computational Mathematics}, 2019
		
		\bibitem{pedrigi}
		Pedrigi R.M., Poulsen C.B., Mehta VV, Ramsing Holm N., Pareek N., Post A.L., Kilic I.D.,
		Banya W.A.S., Dall’Ara G., Mattesini A., Bjørklund M.M., Andersen N.P., Grøndal A.K.,
		Petretto E., Foin N., Davies J.E., Mario C., Di Fog Bentzon J., Erik Bøtker H., Falk E.,
		Krams R., de Silva R.
		Inducing persistent flow disturbances accelerates atherogenesis and promotes thin cap fibroatheroma development in D374Y-PCSK9 hypercholesterolemic minipigs.{\it Circulation}, 132,1003–1012, 2015,
		DOI: \href{https://doi.org/10.1161/CIRCULATIONAHA.115.016270}{10.1161/CIRCULATIONAHA.115.016270}
		
		\bibitem{pfeiffer15}
		Peiffer V., Sherwin S.J., Weinberg P.D.:
		Computation in the rabbit aorta of a new metric—the transverse wall shear stress—to quantify the multidirectional character of disturbed blood flow. {\it Journal of Biomechanics}, 46, 2651–2658, 2013.
		DOI: \href{https://doi.org/10.1016/j.jbiomech.2013.08.003}{10.1016/j.jbiomech.2013.08.003}
		
		\bibitem{pfeiffer} Pfeiffer, P., Sherwin, S.J., Weinberg, P.D.:
		Does low and oscillatory wall shear stress correlate spatially with early atherosclerosis? A systematic review,
		{\it Cardiovascular Research}, 99, 242-250, 2013,
		DOI: \href{https://doi.org/10.1093/cvr/cvt044}{10.1093/cvr/cvt044}
		
		\bibitem{Hundertmark2013} Rusn\'akova G., Luk\'a\v cov\'a, M., Hundertmark, A.:  Kinematic splitting algorithm for fluid-structure interaction in hemodynamics,  {\it Computer Methods in Applied Mechanics and Engineering} 265,  83-106, 2013
		
		
		\bibitem{Giorgio_WSS2}  
		Shojima, M., Oshima, M.,  Takagi, K.,  Torii, R.,  Hayakawa, M.,  Katada, K., Morita, A., Kirino, T.:
		Magnitude and Role of Wall Shear Stress on Cerebral Aneurysm
		Computational Fluid Dynamic Study of 20 Middle Cerebral Artery Aneurysms
		{\it Stroke},  35 (11),  2500–2505,  2004,
		DOI: \href{https://doi.org/10.1161/01.STR.0000144648.89172.0f }{10.1161/01.STR.0000144648.89172.0f}
		
		\bibitem{soulis}
		Soulis, S. V., Giannoglou, G., Fytanidis, D.K.:
		Relative residence time and oscillatory shear index of non-Newtonian flow models in aorta,
		{\it International Workshop on Biomedical Engineering, Biomedical Engineering}, 10, 2013,
		DOI: \href{https://doi.org/10.1109/IWBE.2011.6079011}{10.1109/IWBE.2011.6079011}
		
		
		\bibitem{spanos}
		Spanos, K., Petrocheilou, G., Karathanos, C., Labropoulos, N., Mikhailidis, D., Giannoukas, A.:
		Carotid Bifurcation Geometry and Atherosclerosis 
		{\it Angiology}, 68(9), 757-764, 2017,
		DOI: \href{https://doi-org.wwwdb.dbod.de/10.1177/0003319716678741}{10.1177/0003319716678741}
		
		
		
		
		
		\bibitem{who01} WHO, \textit{The top 10 cause of death},\\ \href{https://www.who.int/news-room/fact-sheets/detail/the-top-10-causes-of-death}{https://www.who.int/news-room/fact-sheets/detail/the-top-10-causes-of-death}, December 2020
		
		
		
		\bibitem{Tremmel_OSI}  Xiang, J., Sabareesh K. Natarajan, S. K., Tremmel, M., Ma, D.,   Mocco, J., Hopkins  L. N.,  Siddiqui, A. H.,
		Levy, E. I., Meng H.,:
		Hemodynamic-Morphologic Discriminants for Intracranial Aneurysm Rupture, {\it Stroke},  42(1),  144–152, 2011, 
		DOI: \href{https://www.ahajournals.org/doi/full/10.1161/STROKEAHA.110.592923}{10.1161/STROKEAHA.110.592923}
		
		\bibitem{zarins_streamlines}
		Zarins, C.K., Giddens, D.P., Bharadvaj, B. K., Sottiurai, V.S.,  Mabon, R.F., Glagov, S.:
		Carotid Bifurcation  Atherosclerosis: Quantitative  Correlation  of Plaque  Localization with  Flow  Velocity Profiles  and Wall  Shear  Stress
		{\it Circulation Research}, 53(4), 502-514, 1983,
		DOI: \href{https://doi.org/10.1161/01.RES.53.4.502}{10.1161/01.RES.53.4.502}
		
		\bibitem{taylor}
		Taylor, C. A., Hughes, T. J. R., Zarins, C. K.: Finite element modeling of three-dimensional pulsatile flow in the abdominal aorta: relevance to atherosclerosis. {\it Annals of Biomedical Engineering}, 26, 975–987, 1998,
		DOI: \href{https://doi.org/10.1114/1.140}{10.1114/1.140}
		
\end{thebibliography}
\end{document}